\documentclass[twoside, 12pt]{article}
\usepackage[latin1]{inputenc}
\usepackage[T1]{fontenc}
\usepackage{amscd}
\usepackage{amsmath}
\usepackage{amsfonts}
\usepackage{amssymb}
\usepackage{amsthm}
\usepackage{comment}

\newtheorem{thm}{Theorem}[section]
\newtheorem*{thm*}{Theorem}

\newtheorem{lem}[thm]{Lemma}

\theoremstyle{definition}

\newtheorem*{rmq}{Remark}

\renewcommand{\tilde}{\widetilde}

\newcommand{\norm}[1]{\left\| #1 \right\|}

    \DeclareMathOperator{\Card}{Card}
    
    \DeclareMathOperator{\diam}{diam}

    \DeclareMathOperator{\Cor}{Cor}
    \DeclareMathOperator{\dist}{dist}
    \DeclareMathOperator{\Leb}{Leb}

\newcommand{\dd}{\, {\rm d}}

\renewcommand{\geq}{\geqslant}
\renewcommand{\leq}{\leqslant}

\newcommand{\N}{\mathbb{N}}

\newcommand{\R}{\mathbb{R}}

\newcommand{\C}{\mathbb{C}}

\renewcommand{\phi}{\varphi}
\renewcommand{\epsilon}{\varepsilon}
\newcommand{\tq}{\ |\ }

\newcommand{\boB}{\mathcal{B}}

\newcommand{\boC}{\mathcal{C}}

\newcommand{\moins}{\backslash}

\newcommand{\boU}{\mathcal{U}}

\title{Decay of correlations for nonuniformly expanding systems
\footnote{\emph{keywords}: decay of correlations, Young tower, non
uniformly expanding maps.
  \emph{2000 Mathematics Subject Classification:} 37A25, 37D25}}

\author{S\'ebastien Gou\"ezel
  \footnote{D\'epartement de Math\'ematiques et Applications,
École Normale Sup\'erieure, 45 rue d'Ulm 75005 Paris (France).
e-mail \texttt{Sebastien.Gouezel@ens.fr}}}

\date{December 2004}
\begin{document}
\maketitle

\begin{abstract}
We estimate the speed of decay of correlations for general nonuniformly
expanding dynamical systems, using estimates on the time the system
takes to become really expanding. Our method can deal with fast
decays, such as exponential or stretched exponential.
 We prove in particular that the
correlations of the Alves-Viana map decay in $O(e^{-c \sqrt{n}})$.
\end{abstract}

\section{Results}

\subsection{Decay of correlations and asymptotic expansion}

When $T:M\to M$ is a map on a compact space, the asymptotic
behavior of Lebesgue-almost every point of $M$ under the iteration
of $T$ is related to the existence of absolutely continuous (or
more generally SRB) invariant probability measures $\mu$. To
understand more precisely the mixing properties of the system, an
essential feature is the speed at which the correlations
  \begin{equation*}
  \Cor(f,g\circ T^n):=\int f\cdot g\circ T^n \dd\mu-\int f\dd\mu
  \int g\dd\mu
  \end{equation*}
tend to $0$. In a uniformly expanding setting, the decay is
exponential, but little is known when the expansion is non uniform.

Recently, \cite{alves_luzzatto_pinheiro} introduced a quantitative
way to measure the non-uniform expansion of a map, and
showed that
this ``measure of non-uniformity'' makes
it possible to control the speed of decay of correlations. More
precisely, when the measure of non-uniformity decays polynomially,
\cite{alves_luzzatto_pinheiro} shows that the decay of
correlations is also polynomial, using hyperbolic times techniques
(\cite{alves:multidim_SRB}) and Young towers (\cite{
lsyoung:recurrence}). As a consequence of this result, the
correlations of the Alves-Viana map (\cite{viana:multidim_attr})
decay faster than any polynomial (which implies for example a
central limit theorem). However, all the estimates of
\cite{viana:multidim_attr} are in $O(e^{-c \sqrt{n}})$, which is
stronger. A precise study of the recurrence makes it in fact
possible to show that the correlations also decay in
$O(e^{-c\sqrt{n}})$ (\cite{baladi_gouezel:viana1},
\cite{baladi_gouezel:viana2}). However, this direct approach
relies strongly on the specificities of the Alves-Viana map,
contrary to the approach of \cite{alves_luzzatto_pinheiro}, which
uses only some general abstract properties, and can therefore be
extended to many other cases. The goal of this article is to
extend the results of \cite{alves_luzzatto_pinheiro} (using a
substantially different method) to speeds of $e^{-c\sqrt{n}}$
(among others), which implies that the results of
\cite{baladi_gouezel:viana2} hold in a much wider setting.

Let $M$ be a compact Riemannian manifold (possibly with boundary)
and $T:M\to M$. We assume that there exists a closed subset
$S\subset M$, with zero Lebesgue measure (containing possibly
discontinuities or critical points of $T$, and with $\partial M
\subset S$), such that $T$ is a $C^2$ local diffeomorphism on
$M\moins S$, and is non uniformly expanding: there exists
$\lambda>0$ such that, for Lebesgue almost every $x\in M$,
  \begin{equation}
  \label{non_unif_hyp}
  \liminf_{n\to \infty} \frac{1}{n}\sum_{k=0}^{n-1} \log
  \norm{DT(T^k x)^{-1}}^{-1} \geq \lambda.
  \end{equation}

We also need non-degeneracy assumptions close to $S$, similar to the
assumptions in \cite{alves_bonatti_viana} or
\cite{alves_luzzatto_pinheiro}:
we assume that there exist $B>1$ and $\beta>0$ such that,
for any $x\in
M\moins S$ and every $v\in T_x M \moins \{0\}$,
  \begin{equation}
  \frac{1}{B}\dist(x,S)^\beta \leq \frac{\norm{DT(x)v}}{\norm{v}}
  \leq B \dist(x,S)^{-\beta}.
  \end{equation}
Assume also that, for all $x,y\in M$ with
$\dist(x,y)<\dist(x,S)/2$,
  \begin{equation}
  \Bigl|\log \norm{DT(x)^{-1}} -\log \norm{DT(y)^{-1}} \Bigr|
  \leq B \frac{\dist(x,y)}{\dist(x,S)^\beta}
  \end{equation}
and
  \begin{equation}
  \bigl|\log |\det DT(x)^{-1}| - \log |\det DT(y)^{-1}| \bigr|
  \leq B \frac{\dist(x,y)}{\dist(x,S)^\beta},
  \end{equation}
i.e.\ $\log \norm{DT^{-1}}$ and $\log |\det DT^{-1}|$ are locally
Lipschitz, with a constant which is controlled  by the distance to the
critical set. This implies that the singularities are at most
polynomial, and in particular that the critical points are not flat.

We assume that the critical points come subexponentially close to $S$
in the following sense. For  $\delta>0$, set
$\dist_{\delta}(x,S)=\dist(x,S)$ if $\dist(x,S)<\delta$, and
$\dist_{\delta}(x,S)=1$ otherwise. We assume that, for all
$\epsilon>0$, there exists $\delta(\epsilon)>0$ such that, for
Lebesgue almost every $x\in M$,
  \begin{equation}
  \label{conv_lente}
  \limsup_{n\to \infty} \frac{1}{n}\sum_{k=0}^{n-1}- \log
  \dist_{\delta(\epsilon)}(T^k x,S) \leq \epsilon.
  \end{equation}

We will need to control more precisely the speed of convergence in
\eqref{non_unif_hyp} and \eqref{conv_lente}. As
\cite{alves_luzzatto_pinheiro}, we consider for this the following
function, which measures the non-uniformity of the system:
  \begin{multline*}
  h^1_{(\epsilon_1,\epsilon_2)}(x)
  =\inf \Biggl\{ N \in \N^* \tq \forall n \geq N,
  \frac{1}{n}\sum_{k=0}^{n-1} \log
  \norm{DT(T^k x)^{-1}}^{-1} \geq \frac{\lambda}{2}
  \\
  \text{ and for }i=1,2,
  \frac{1}{n}\sum_{k=0}^{n-1}
  - \log \dist_{\delta(\epsilon_i)}(T^k x,S) \leq 2\epsilon_i \Biggr\}.
  \end{multline*}
It is important to have two indexes $\epsilon_1$ and $\epsilon_2$
to guarantee the existence of hyperbolic times (see Lemma
\ref{existe_epsilon_1_2}). To simplify the notations, we will
write $\epsilon=(\epsilon_1,\epsilon_2)$. The points $x$ such that
$h^1_\epsilon(x)=n$ are ``good'' for times larger than $n$. Hence,
the lack of expansion of the system at time $n$ is evaluated by
  \begin{equation}
  \label{mauvais_points}
  \Leb \{ x \tq h^1_\epsilon(x) > n\},
  \end{equation}
and it is natural to try to estimate the speed of decay of
correlations using this quantity. This is done in
 \cite{alves_luzzatto_pinheiro} in the polynomial case:
if $\eqref{mauvais_points}=O(1/n^\gamma)$ for some
$\gamma>1$, then the correlations of H\"older functions decay at least
like
$1/n^{\gamma-1}$.

Set $\Lambda =\bigcap_{n\geq 0} T^n(M)$. We will say that $T$ is
topologically transitive on the attractor $\Lambda$ if, for every
nonempty open subsets $U,V$ of $\Lambda$, there exists $n$ such
that $T^{-n}(U) \cap V$ contains a nonempty open set (the precise
formulation is important since $T$ may not be continuous on $S$).

We will say that a sequence  $(u_n)_{n\in \N}$ has  \emph{polynomial
decay} if there exists $C>0$ such that, for all $n/2 \leq k \leq
n$, $0< u_k \leq C u_n$. This implies in particular that $u_n$ does
not tend too fast to $0$: there exists $\gamma>0$ such that
$1/n^\gamma=O(u_n)$ (for example $\gamma=\frac{\log C}{\log 2}$).

Finally, the \emph{basin} of a probability measure $\mu$ is the set of
points $x$ such that $\frac{1}{n}\sum_{k=0}^{n-1}\delta_{T^k x}$
converges weakly to $\mu$, where $\delta_y$ is the Dirac mass at $y$.

\begin{thm}
\label{thm_principal_0}
Assume that all the iterates of $T$ are topologically transitive on
$\Lambda$ and that, for all $\epsilon=(\epsilon_1,\epsilon_2)$, there
exists a sequence $u_n(\epsilon)$ with $\sum u_n(\epsilon)<+\infty$
and
  $
  \Leb \{ x \tq h^1_\epsilon(x) > n\}
  =O(u_n(\epsilon)).
  $
Assume moreover that $u_n(\epsilon)$ satisfies one of the following properties:
\begin{enumerate}
\item
\label{cas_1_kjsaf}
$u_n(\epsilon)$ has polynomial decay.
\item
\label{cas_2_kjsaf}
There exist $c(\epsilon)>0$ and $\eta(\epsilon) \in (0,1]$ such that
$u_n(\epsilon)=e^{-c(\epsilon) n^{\eta(\epsilon)}}$.
\end{enumerate}
Then $T$ preserves a unique (up to normalization) absolutely
continuous (with respect to Lebesgue) measure $\mu$. Moreover, this
is a mixing probability measure, whose basin contains Lebesgue-almost
every point of $M$.

Finally, there exists $\epsilon^0=(\epsilon^0_1,\epsilon^0_2)$
such that, if $f,g:M \to \R$ are
two functions with $f$ H\"older and $g$ bounded,
their correlations
$\Cor(f,g\circ
T^n)=\int f\cdot g\circ T^n \dd\mu -\int f \dd \mu \int g \dd \mu$
decay at the following speed:
\begin{enumerate}
\item $|\Cor(f,g\circ T^n)| \leq C \sum_{p=n}^\infty u_p(\epsilon^0)$
in case \ref{cas_1_kjsaf}.
\item There exists $c'>0$ such that
$|\Cor(f,g\circ T^n)| \leq C e^{-c' n^{\eta(\epsilon^0)}}$ in case
 \ref{cas_2_kjsaf}.
\end{enumerate}
\end{thm}
In fact, $\epsilon^0$ can be chosen \emph{a priori}, depending
only on  $\lambda$ and $T$. It would then be sufficient to have
\eqref{conv_lente} for $\epsilon^0_1$ and $\epsilon^0_2$ to get
the theorem. However, in practical cases, it is often not harder
to prove \eqref{conv_lente} for all values of $\epsilon$ than to
prove it for a specific value of $\epsilon$. This is why, as in
\cite{alves_bonatti_viana} and \cite{alves_luzzatto_pinheiro}, we
have preferred to state the theorem in this more convenient way.

In the first case, taking $u_n=1/n^\gamma$, we get another proof of the result
of \cite{alves_luzzatto_pinheiro}.

The main problem of this theorem is that \eqref{mauvais_points} is
often difficult to estimate, since $h^1_\epsilon(x)$ states a
condition on \emph{all} iterates of $x$, and not only a finite number
of them.

\subsection{The Alves-Viana map}

Theorem \ref{thm_principal_0} applies to the Alves-Viana map, given by
  \begin{equation}
  T : \left\{ \begin{array}{ccc}
  S^1 \times I & \to & S^1 \times I \\
  (\omega,x) & \mapsto & (16 \omega, a_0+\epsilon \sin(2\pi
  \omega)-x^2)
  \end{array} \right.
  \end{equation}
where $a_0 \in (1,2)$ is a Misiurewicz point (i.e.\ the critical point
$0$ is preperiodic for  $x\mapsto
a_0-x^2$), $\epsilon$ is small enough and $I$ is a compact interval of
$(-2,2)$ such that $T$ sends $S^1\times I$ into its interior.

This map has been introduced by Viana in
\cite{viana:multidim_attr}. He shows that $T$ (and in fact any map
close enough to $T$ in the $C^3$ topology) has almost everywhere
two positive Lyapunov exponents, even though there are critical
points in the fibers. More precisely, Viana shows that the points
that do not see the expansion in the fiber have a measure decaying
like $O(e^{-c \sqrt{n}})$.
\cite{alves_araujo:random_perturbations} obtains from this
information that, for every $\epsilon=(\epsilon_1,\epsilon_2)$,
for every $c<1/4$,
  \begin{equation}
  \label{eq_viana_O}
  \Leb\{ x \tq h^1_\epsilon(x) > n\}=O(e^{-c \sqrt{n}}).
  \end{equation}
Moreover, \cite{alves_viana} shows that all the iterates of $T$ are
topologically transitive on $\Lambda$.

A consequence of the results of \cite{alves_luzzatto_pinheiro} is
that the correlations of the Alves-Viana map decay faster than any
polynomial. However, their method of proof can deal only with
polynomial speeds (see paragraph \ref{par_1.4_jasklf}), and hence
can not reach the conjectural upper bound of $e^{-c' \sqrt{n}}$.
Theorem \ref{thm_principal_0} implies this conjecture (already
announced in \cite{baladi_gouezel:viana1}):

\begin{thm}
\label{decor_viana_rapide}
The correlations of H\"older functions for any map close
enough (in the $C^3$ topology) to the Alves-Viana map decay at least
like $e^{-c'\sqrt{n}}$ for some $c'>0$.
\end{thm}
This result applies also if the expansion coefficient $16$ is
replaced by $2$, according to \cite{buzzi_sester_tsujii}.
Note that the specific method of \cite{baladi_gouezel:viana2}, which
proves Theorem \ref{decor_viana_rapide}, can not be directly used when
$16$ is replaced by $2$, since it uses in particular the specific form
of admissible curves. On the other hand, the abstract method of this
article applies immediately, since  \cite{buzzi_sester_tsujii} proves
essentially \eqref{eq_viana_O}.

\subsection{Decorrelation and expansion in finite time}

The function $h^1_\epsilon(x)$ takes into account the expansion at $x$
for large enough times, and is consequently hard to estimate in
general. It is more natural to consider the first time with enough
expansion. For technical reasons, we will need three parameters to get
results in this setting (see the proof of Lemma
\ref{existe_epsilon_3}). Set
  \begin{multline*}
  h^2_{(\epsilon_1,\epsilon_2,\epsilon_3)}(x)
  =\inf \Biggl\{ n \in \N^* \tq
  \frac{1}{n}\sum_{k=0}^{n-1} \log
  \norm{DT(T^k x)^{-1}}^{-1} \geq \frac{\lambda}{2}\\
  \text{ and for }i=1,2,3,
  \frac{1}{n}\sum_{k=0}^{n-1}
  - \log \dist_{\delta(\epsilon_i)}(T^k x,S) \leq 2\epsilon_i \Biggr\}.
  \end{multline*}
This definition takes only the first $n$ iterates of $x$ into account,
and can consequently be checked in finite time. We will write
$\epsilon=(\epsilon_1,\epsilon_2,\epsilon_3)$. The time
$h^2_\epsilon$ is related to the notion of \emph{first hyperbolic
time} studied for example in  \cite{alves_araujo:hyperbolic_times}.

If there were only two parameters in the definition of $h^2$, we would
have $h^2\leq h^1$. However, since there are three parameters, $h^1$
and $h^2$ can rigourously not be compared.

We will estimate the speed of decay of correlations using $\Leb\{x
\tq h^2_\epsilon(x)>n\}$. Our main result is the following
theorem:
\begin{thm}
\label{thm_principal}
Assume that all the iterates of $T$ are topologically transitive on
$\Lambda$ and that, for all
$\epsilon=(\epsilon_1,\epsilon_2,\epsilon_3)$, there exists a sequence
$u_n(\epsilon)$ with $\sum (\log n)u_n(\epsilon)<+\infty$ and
  $
  \Leb\left\{ x \tq h^2_\epsilon(x)>n\right\}
  =O(u_n(\epsilon)).
  $
Assume moreover that $u_n(\epsilon)$ satisfies one of the following properties:
\begin{enumerate}
\item
\label{laksjdflk;jasdf}
$u_n(\epsilon)$ has polynomial decay.
\item
\label{aksfmnavs,sadf}
there exist $c(\epsilon)>0$ and $\eta(\epsilon) \in (0,1]$ such that
$u_n(\epsilon)=e^{-c(\epsilon) n^{\eta(\epsilon)}}$.
\end{enumerate}
Then $T$ preserves a unique (up to normalization) absolutely
continuous invariant measure $\mu$. Moreover, this measure is a mixing
probability measure, whose basin contains Lebesgue almost every point
of $M$.

Finally, there exists
$\epsilon^0= (\epsilon^0_1, \epsilon^0_2,
\epsilon^0_3)$ such that,
if $f,g:M \to \R$ are
two functions with $f$ H\"older and $g$ bounded,
their correlations
$\Cor(f,g\circ T^n)=\int f\cdot
g\circ T^n \dd\mu -\int f \dd \mu \int g \dd \mu$ decay at the
following speed:
\begin{enumerate}
\item $|\Cor(f,g\circ T^n)| \leq C \sum_{p=n}^\infty (\log
p)u_p(\epsilon^0)$ in case \ref{laksjdflk;jasdf}.
\item There exists $c'>0$ such that
$|\Cor(f,g\circ T^n)| \leq C e^{-c' n^{\eta(\epsilon^0)}}$ in case
\ref{aksfmnavs,sadf}.
\end{enumerate}
\end{thm}
For example, when
$\Leb\{ x \tq h^2_\epsilon(x)>n \} =
O(1/n^\gamma)$ with $\gamma>1$, the correlations decay like
$\log n/n^{\gamma-1}$. In the first case (polynomial decay), note
that there is a loss of $\log n$ between
Theorem \ref{thm_principal_0} and Theorem \ref{thm_principal}. It is
not clear whether this loss is real, or due
to the technique of proof.

The comments on the choice of $\epsilon^0$ following Theorem
\ref{thm_principal_0} are still valid here. It is even possible to
take the same value for $\epsilon^0_1$ and $\epsilon^0_2$ in both theorems.

We will return later to the existence of invariant measures
(Theorems \ref{description_mesure_invariante} and
\ref{vitesse_decroissance}). Without transitivity assumptions, we will
get a spectral decomposition: $T$ admits a finite number of absolutely
continuous invariant ergodic probability measures, and each of these
measures has a finite number of components which are mixing for an
iterate of $T$, with the same bounds on the decay of correlations as
in Theorems \ref{thm_principal_0} and \ref{thm_principal}: these
theorems correspond to the case where the spectral decomposition is trivial.

\begin{rmq}
If $u_n$ has polynomial decay and
$u_n=O(1/n^\gamma)$
for some $\gamma>1$, then $\sum_{p=n}^\infty (\log
p)u_p=O\left((\log n)\sum_{p=n}^\infty u_p\right)$, which simplifies a
little the bound on the decay of correlations.
\end{rmq}

\begin{rmq}
In the stretched exponential case (i.e.\ $0<\eta<1$), the conclusions
of Theorems \ref{thm_principal_0} and
\ref{thm_principal} are true for any $c'<c(\epsilon^0)$. This can easily
be checked in all the following proofs (except in the proof of Lemma
\ref{lemme_ameliore_Young}, where slightly more careful estimates are
required).
\end{rmq}

\subsection{Strategy of proof}

\label{par_1.4_jasklf}
As it is often the case when one wants to estimate the decay of
correlations, the strategy of proof consists in building a Young
tower (\cite{lsyoung:recurrence}), i.e.\ selecting a subset $B$ of
$M$ and building a partition $B=\bigcup B_i$ such that $T^{R_i}$
is an isomorphism between $B_i$ and $B$, for some return time
$R_i$. Then \cite{lsyoung:recurrence} gives estimates on the decay
of correlations, depending on the measure of points coming back to
$B$ after time $n$, i.e., $\Leb\left(\bigcup_{R_i>n} B_i \right)$.
To construct the sets $B_i$, we will use \emph{hyperbolic times}.
Denote by $H_n$ the set of points for which $n$ is a hyperbolic
time.

This strategy is implemented in \cite{alves_luzzatto_pinheiro}. We
will describe quickly their inductive construction, in a somewhat
incorrect way but giving the essential ideas. Before time $n$,
assume that some sets $B_i$ have already been constructed, with a
return time $R_i$ satisfying $R_i<n$. At time $n$, consider $H_n
\moins \left(\bigcup_{R_i<n} B_i\right)$, and construct new sets
$B_j$ covering a definite proportion of this set, with return time
$R_j=n$. Using some information about the repartition of
hyperbolic times (the \emph{Pliss Lemma}), it is then possible to
prove that $\Leb\left(\bigcup_{R_i>n} B_i \right)$ decays at least
polynomially. The main limitation of this strategy is that, at
time $n$, it can deal only with a fraction of $H_n$. Since the
repartition of hyperbolic times is \emph{a priori} unknown (except
for the Pliss Lemma), we may have to wait a long time ($\sim n$)
to see another hyperbolic time. This makes it impossible to prove
that the decorrelations decay faster than $e^{-c (\log n)^2}$
without further information.

To avoid this problem, we will deal with all points of $H_n$ at
time $n$, and not only a fraction. To do this, we will consider a
fixed partition $U_1,\ldots, U_N$ of the space (with $N$ fixed)
and use $f^{-n}(U_1),\dots, f^{-n}(U_N)$ to partition $H_n$. In
this way, we will get a partition $\boB_i$ of $U_i$ (for each
$i$), and each element of $\boB_i$ will be sent on some (possibly
different) $U_j$ by an iterate of $T$. Moreover, we will keep a
precise control on the measure of points having long return times.

Using this auxiliary partition, it will be quite easy to build a
Young tower, using an inductive process: select some $U_i$, for
example $U_1$. While a point does not fall into $U_1$, go on
iterating, so that it falls in some $U_j$, then some $U_k$, and so
on. Most points will come back to $U_1$ after a finite (and well
controlled) number of iterates, and this will give the required
partition of $U_1$.

Finally, to estimate the decay of correlations, it will not be
possible to
apply directly the results of \cite{lsyoung:recurrence}, since they are
slightly too weak
(in the case of $e^{-c n^\eta}$ with $0<\eta<1$, Young proves
only a decay of correlations of
$e^{-c' n^{\eta'}}$ for any $\eta'<\eta$, which is weaker than the
results of Theorems \ref{thm_principal_0} and \ref{thm_principal}).
However, the combinatorial techniques used in the construction of the
partition will easily enable us to strengthen the results of
\cite{lsyoung:recurrence}, to obtain the required estimates.

The main difficulty of the proof will be to get the estimates on the
auxiliary partition $U_1,\ldots,U_N$, in Section \ref{section_Ui} (for
example, the logarithmic loss between Theorems \ref{thm_principal} and
\ref{thm_principal_0} will appear there). Then we will build the
Young tower in Section \ref{section_tour}, and estimate the decay of
correlations in paragraph \ref{section_decorrelation}. We will prove at
the same time Theorems \ref{thm_principal_0} and \ref{thm_principal}.

\textbf{Acknowledgments.} I would like to thank V. Baladi for many
enlightening discussions and explanations, and the referee for his
valuable comments.

\newpage

\section{Hyperbolic times}
\label{section_temps_hyperboliques}

We recall in this section the notion of hyperbolic times, of
\cite{alves:multidim_SRB} and \cite{alves_bonatti_viana}, and we
describe different sets that can be built at hyperbolic times.
These sets will be the basic stones used to build the auxiliary partition
in Section
\ref{section_Ui}.

Let $b$ be a constant such that $0<b<\min(1/2,1/(4\beta))$.
For $\sigma<1$ and $\delta>0$, we say that $n$ is a
$(\sigma,\delta)$-hyperbolic time for $x$ if, for all $1\leq k\leq n$,
  \begin{equation}
  \prod_{j=n-k}^{n-1}\norm{DT(T^j x)^{-1}} \leq \sigma^k
  \text{ and } \dist_\delta(T^{n-k}x,S) \geq \sigma^{b k}.
  \end{equation}
We will denote by
$H_n=H_n(\sigma,\delta)$ the set of points for which $n$ is a
$(\sigma, \delta)$-hyperbolic time.

In paragraph \ref{section_frequence_hyperboliques}, we will choose
carefully the constants $\sigma$ and $\delta$ (as well as $\epsilon^0$
given by Theorems \ref{thm_principal_0} and
\ref{thm_principal}). However, the reasons for this choice will not
become clear before paragraph \ref{prouve_thm_auxiliaire}, and the
reader may admit the existence of $\sigma, \delta$ and
$\epsilon_0$, and come back to paragraph
\ref{section_frequence_hyperboliques} just before reading paragraph
\ref{prouve_thm_auxiliaire}.

\subsection{Frequency of hyperbolic times}

\label{section_frequence_hyperboliques}

The following lemma is a slight generalization of
\cite[Lemma
5.4]{alves_bonatti_viana}:

\begin{lem}
\label{existe_epsilon_3}
Take $T:M\to M$ and $\delta:\R_+^*\to \R_+^*$ such that
\eqref{non_unif_hyp} and \eqref{conv_lente} are satisfied.

Then there exist $\epsilon_3>0$ and $\kappa>0$ such that, for all
$\epsilon_1,\epsilon_2< \kappa$, there exists
$\theta(\epsilon_1,\epsilon_2)>0$ such that, if $x\in M$ and $n\in
\N^*$ satisfy
  \begin{equation*}
  \frac{1}{n}\sum_{k=0}^{n-1} \log
  \norm{DT(T^k x)^{-1}}^{-1} \geq \frac{\lambda}{2}
  \text{ and for }i=1,2,3,
  \frac{1}{n}\sum_{k=0}^{n-1}
  - \log \dist_{\delta(\epsilon_i)}(T^k x,S) \leq 2\epsilon_i,
  \end{equation*}
then there exist times $1\leq p_1<\ldots<p_l \leq n$ with
$l\geq \theta(\epsilon_1,\epsilon_2) n$ such that, for all $j\leq l$,
  \begin{multline}
  \label{super_temps_hyperb}
  \forall 1 \leq k\leq p_j,
  \sum_{s=p_j-k}^{p_j-1} \log
  \norm{DT(T^s x)^{-1}}^{-1} \geq \frac{\lambda}{4}k\\
  \text{ and for }i=1,2,
  \sum_{s=p_j-k}^{p_j-1}
  - \log \dist_{\delta(\epsilon_i)}(T^s x,S) \leq 2\sqrt{\epsilon_i} k.
  \end{multline}
\end{lem}
This means that the density of times $p$ between $1$ and $n$
satisfying \eqref{super_temps_hyperb} is at least
$\theta(\epsilon_1,\epsilon_2)$.
Before giving the proof of the lemma, we will state another lemma with
the same flavor:
\begin{lem}
\label{existe_epsilon_1_2}
Take $T:M\to M$ and $\delta:\R_+^*\to \R_+^*$ such that
\eqref{non_unif_hyp} and \eqref{conv_lente} are satisfied. Take also
$\kappa>0$.

Then there exist $\epsilon_1, \epsilon_2 < \kappa$ and $\theta>0$ such
that,
if $x \in M$ and $n\in \N^*$ satisfy
  \begin{equation*}
  \frac{1}{n}\sum_{k=0}^{n-1} \log
  \norm{DT(T^k x)^{-1}}^{-1} \geq \frac{\lambda}{4}
  \text{ and for }i=1,2,
  \frac{1}{n}\sum_{k=0}^{n-1}
  - \log \dist_{\delta(\epsilon_i)}(T^k x,S) \leq 2\sqrt{\epsilon_i},
  \end{equation*}
then there exist times $1\leq p_1<\ldots<p_l \leq n$ with
$l\geq \theta n$ such that, for all $j\leq l$,
  \begin{equation*}
  \forall 1 \leq k\leq p_j,
  \sum_{s=p_j-k}^{p_j-1} \log
  \norm{DT(T^s x)^{-1}}^{-1} \geq \frac{\lambda}{8}k
  \text{ and }
  \sum_{s=p_j-k}^{p_j-1}
  - \log \dist_{\delta(\epsilon_1)}(T^s x,S) \leq b \frac{\lambda}{8} k.
  \end{equation*}
\end{lem}

Until the end of this article, we will denote by $\epsilon_3^0$
the value of $\epsilon_3$ given by Lemma \ref{existe_epsilon_3},
and by $\epsilon_1^0,\epsilon_2^0$ the values of $\epsilon_1$ and
$\epsilon_2$ given by Lemma \ref{existe_epsilon_1_2}. We will also
set $\sigma=e^{-\lambda/8}<1$. Finally, write
$\delta=\delta(\epsilon_1^0)$.

Hence, the times $p_j$ given by the conclusion of Lemma
\ref{existe_epsilon_1_2} are
$(\sigma,\delta)$-hyperbolic. In the same way, the times $p_j$
satisfying \eqref{super_temps_hyperb} are also
$(\sigma,\delta)$-hyperbolic (if $\kappa$ is small enough), but they
are more than that since they guarantee a control at the same time for
 $\epsilon_1^0$ and for $\epsilon_2^0$ (whence Lemma
\ref{existe_epsilon_1_2} can be applied to them): we will say that a
time which satisfies
\eqref{super_temps_hyperb} for $\epsilon^0_1$ and $\epsilon^0_2$
is a \emph{super hyperbolic time}.
We will write $SH_n$  for the set of points for which $n$ is a super
hyperbolic time, and $H_n=H_n(\sigma,\delta)$ for the set of points
for which $n$ is a  $(\sigma,\delta)$-hyperbolic time. In particular,
$SH_n \subset H_n$.

In the following proof, we will see why an index $\epsilon$ is lost:
it is used to obtain the conclusion on $\sum_{s=p_j-k}^{p_j-1} \log
  \norm{DT(T^s x)^{-1}}^{-1}$, since Pliss Lemma can not be applied
directly (since this sequence is not bounded), whence another control
is needed.

\begin{proof}[Proof of Lemma \ref{existe_epsilon_3}]
The proof is essentially the proof of Lemma 5.4 of
\cite{alves_bonatti_viana}: they first show that there exist
$\epsilon_3>0$ (which can be taken arbitrarily small) and $\theta_1>0$
such that, if
  \begin{equation*}
  \frac{1}{n}\sum_{k=0}^{n-1} \log
  \norm{DT(T^k x)^{-1}}^{-1} \geq \frac{\lambda}{2}
  \text{ and }
  \frac{1}{n}\sum_{k=0}^{n-1}
  - \log \dist_{\delta(\epsilon_3)}(T^k x,S) \leq 2\epsilon_3,
  \end{equation*}
then there is a proportion at least $\theta_1>0$ of times $p$ between
$1$ and $n$ such that
  \begin{equation*}
  \forall 1 \leq k\leq p,
  \sum_{s=p-k}^{p-1} \log
  \norm{DT(T^s x)^{-1}}^{-1} \geq \frac{\lambda}{4}k.
  \end{equation*}
Moreover, \cite[Lemma 3.1]{alves_bonatti_viana} also shows that,
for $\epsilon>0$, if $x$ satisfies
  \begin{equation*}
  \frac{1}{n}\sum_{k=0}^{n-1}
  - \log \dist_{\delta(\epsilon)}(T^k x,S) \leq 2\epsilon,
  \end{equation*}
then there exists a proportion at least
$\theta(\epsilon)=1-\sqrt{\epsilon}$ of times $p$ between $1$ and $n$
such that
   \begin{equation*}
  \forall 1 \leq k\leq p,
  \sum_{s=p-k}^{p-1}
  - \log \dist_{\delta(\epsilon)}(T^k x,S) \leq 2\sqrt{\epsilon}k.
  \end{equation*}
When $\epsilon \to 0$, $\theta(\epsilon) \to 1$. Hence, if $\kappa$ is
small enough, for all
$\epsilon_1,\epsilon_2<\kappa$, we will have
$\theta(\epsilon_1,\epsilon_2):= \theta_1 + \theta(\epsilon_1)
+\theta(\epsilon_2) -2 >0$, which gives the conclusion of the lemma.
\end{proof}
The proof of Lemma \ref{existe_epsilon_1_2} is similar.

\subsection{Constructions at hyperbolic times}

The following lemma refines \cite[Lemma 5.2]{alves_bonatti_viana}
and \cite[Lemma 4.1]{alves_luzzatto_pinheiro}:

\begin{lem}
\label{prop_temps_hyperb}
There exist $\delta_2,D_1,\lambda_1<1$ such that, if $x\in M$ and
$n$ is a $(\sigma,\delta)$-hyperbolic time for $x$, there exists a
unique neighborhood $V_n(x)$ of $x$ with the following properties:
\begin{enumerate}
\item $T^n$ is a diffeomorphism between $V_n(x)$ and the ball
$B(T^nx,\delta_2)$.
\item For $1\leq k\leq n$ and $y,z\in V_n(x)$, $\dist(T^{n-k}y,T^{n-k}z)
\leq \sigma^{k/2} \dist(T^n y,T^n z)$.
\item For all $y,z\in V_n(x)$,
  \begin{equation*}
  \left|\frac{\det DT^n(y)}{\det DT^n(z)}-1 \right| \leq D_1 \dist(T^n
  y, T^n z).
  \end{equation*}
\item $V_n(x) \subset B(x,\lambda_1^n)$.
\item If $n\leq m$, $y\in H_m$ and
$V_n(x) \cap V_m(y) \not=\emptyset$, then $T^n$ is injective on
$V_n(x) \cup V_m(y)$.
\end{enumerate}
\end{lem}
Note that the third assertion of the lemma implies that the
volume-distortion of $T^n$ is bounded by $
D_2:=2\delta_2 D_1+1$, i.e., for all
$U,V\subset V_n(x)$,
  \begin{equation}
  \label{distortion_bornee}
  D_2^{-1} \frac{\Leb (T^n(U))}{\Leb (T^n(V))} \leq \frac{\Leb (U)}{\Leb (V)}
  \leq D_2 \frac{\Leb (T^n(U))}{\Leb (T^n(V))}.
  \end{equation}
\begin{proof}
Lemma 5.2 of \cite{alves_bonatti_viana} shows that there exists
$\delta_1>0$ such that, if $x\in H_n(\sigma, \delta)$, then there
exists a neighborhood
$V'_n(x)$ mapped diffeomorphically by $T^n$ to
$B(T^n x, \delta_1)$. We set $V_n(x)=V'_n(x) \cap
T^{-n} (B(T^n x, \delta_1/4))$, and $\delta_2=\delta_1/4$. As
$V_n(x) \subset V'_n(x)$, the first and second assertion of the lemma
come from Lemma 5.2 of \cite{alves_bonatti_viana}, and the third one
from Lemma 4.1 of \cite{alves_luzzatto_pinheiro}. The fourth one is a
consequence of the second one (for
$\lambda_1=\sigma^{1/2}$).

For the uniqueness, note that two distinct neighborhoods
$V_n^1(x)$ and $V_n^2(x)$ would give two different lifts by $T^n$
of a path from $T^n(x)$ to a point in $B(T^n(x),\delta_1/4)$,
which is not possible.

Finally, assume that $V_n(x) \cap V_m(y)$ contains a point
$z$. Then
  \begin{equation*}
  \diam(T^n(V_m(y))) \leq \diam(T^m(V_m(y))) =\delta_1/2,
  \end{equation*}
whence $T^n(V_m(y)) \subset B(T^n x, \delta_1)$. We build a set
$W_m(y)=T^{-n}( T^n(V_m(y))) \cap V'_n(x)$: by definition of
$V'_n(x)$, $T^n$ is an isomorphism between $W_m(y)$ and
$T^n(V_m(y))$. But $T^n$ is also an isomorphism between  $V_m(y)$
and $T^n(V_m(y))$. As $V_m(y)$ and $W_m(y)$ both contain $z$, the
previous uniqueness argument implies that $V_m(y)=W_m(y)$. In
particular, $V_m(y) \subset V'_n(x)$. As $T^n$ is injective on
$V'_n(x)$, it is also injective on
 $V_n(x) \cup V_m(y)$.
\end{proof}

Take $\boU=\{U_1,\ldots,U_N\}$ a finite partition of $M$ by sets of
diameter at most $\delta_2/10$, with nonempty interiors and piecewise smooth
boundaries (for example a triangulation of $M$). Hence, there
exist constants $C_2>0$ and $\lambda_2<1$ such that
  \begin{equation}
  \label{frontiere_exp}
  \forall 1\leq i \leq N, \forall n\in \N, \Leb\{ x\in U_i \tq
  \dist(x,\partial U_i) \leq \lambda_1^n \} \leq C_2 \lambda_2^n.
  \end{equation}
We will write $U'_i=\{ x\in M \tq \dist(x, U_i) \leq
\delta_2/10\}$. Increasing $C_2$ and $\lambda_2$ if necessary, we can
also assume that
  \begin{equation}
  \label{mesure_frontiere_U'_i}
  \forall 1\leq i \leq N, \forall n\in \N, \Leb\left\{ x\in M \tq
  \dist(x,\partial U'_i) \leq \frac{\delta_2}{2}
  \sigma^{n/2} \right\} \leq C_2 \lambda_2^n \Leb(U_i).
  \end{equation}
We will finally assume that, for any ball $B(x,\delta_2)$ of radius
$\delta_2$ and for all $1\leq i \leq N$,
  \begin{equation}
  \label{grosse_mesure_Ui}
  \Leb B(x,\delta_2) \leq
  C_2 \Leb(U_i).
  \end{equation}

Take $x\in H_n$. Then $T^n x$ belongs to a unique $U_i=:U(x,n)$,
included in $B(T^n x, \delta_2)=T^n(V_n(x))$. We will write
$I^n_\infty(x)=T^{-n}(U_i) \cap V_n(x)$. In the construction of the
auxiliary partition in Section \ref{section_Ui}, the partition
elements will be such sets $I^n_\infty(x)$. In the construction, if we choose
$I^n_\infty(x)$ and then $I^{n+1}_\infty(y)$ while
$y\not\in I^n_\infty(x)$ but $y$ is very close to the boundary of
$I^n_\infty(x)$, the two sets $I^n_\infty(x)$ and
$I^{n+1}_\infty(y)$ may have a nonempty intersection, which we want to
avoid since we are building a partition. As in  \cite{lsyoung:annals},
we will have to introduce a waiting time telling when it is not
dangerous to select $y$, ensuring that
$I^n_\infty(x)
\cap I^m_\infty(y)=\emptyset$. We thus set, for $m>n$,
  \begin{equation*}
  I^n_m(x)=
  \left\{ y\in V_n(x) \tq \frac{\delta_2}{10} \sigma^{\frac{m-n}{2}} <
  \dist(T^n y, U(x,n)) \leq \frac{\delta_2}{10} \sigma^{\frac{m-n-1}{2}}
  \right\}
  \end{equation*}
and $I^n_{\geq m}(x)=\bigcup_{m\leq t<\infty} I^n_t(x)$: these are
the points which are not allowed to be selected at time $m$,
because they could interfere with $x$ at time $n$ (this choice
will be justified by Lemma \ref{lemme_disjointitude}, and
\eqref{eq_disjointitude}). We will say that a point of $I^n_{\geq
m}(x)$ is forbidden by the time $n$, at the time $m$. We will also
write  $\tilde{I}^n_{\geq m}(x)=\bigcup_{m\leq t\leq \infty}
I^n_t(x)$, i.e.\ we add the ``core'' $I^n_\infty(x)$. The main
difference with \cite{lsyoung:recurrence} or
\cite{alves_luzzatto_pinheiro} is that, in these articles, the
combinatorial estimates are less precise, whence they can afford
to forget the time by which a point is forbidden (the $n$ in
$I^n_{\geq m}$).

\begin{lem}
\label{tout_dans_Vn}
If $0< n\leq m$ and $\tilde{I}^n_{\geq n+1}(x) \cap
\tilde{I}^m_{\geq m+1}(y) \not=\emptyset$, then $
\tilde{I}^n_{\geq n+1}(x) \cup \tilde{I}^m_{\geq m+1}(y) \subset
V_n(x)$.
\end{lem}
Note that, when we write $\tilde{I}^n_{\geq n+1}(x)$ (for example), it
is implicit that this set is well defined, i.e.\ that
$x\in H_n$.
\begin{proof}
Take $z\in \tilde{I}^n_{\geq n+1}(x) \cap \tilde{I}^m_{\geq
m+1}(y)$. By Lemma \ref{prop_temps_hyperb},
$T^n( \tilde{I}^m_{\geq m+1}(y)) \subset B(T^n z, \delta_2/2)
\subset B(T^n x, \delta_2)$. In particular, every $u\in
T^n(\tilde{I}^m_{\geq m+1}(y))$ has a preimage $u'$ under $T^n$ in
$V_n(x)$. We have to see that $u'$ belongs to $\tilde{I}^m_{\geq
m+1}(y)$. Otherwise, $u$ would have another preimage $u''$ in
$\tilde{I}^m_{\geq m+1}(y)$. As $V_n(x)\cap V_m(y)$ contains $z$,
the fifth assertion of Lemma \ref{prop_temps_hyperb} gives that
$T^n$ is injective on $V_n(x) \cup V_m(y)$. This is a
contradiction since $u'\not=u''$ but $T^n(u')=T^n(u'')$.
\end{proof}

\begin{lem}
\label{lemme_disjointitude}
There exists $P>0$ such that, for $0< n<m$, $x\in H_n$ and $y\in H_m
\moins \tilde{I}^n_{\geq m}(x)$,
  \begin{equation*}
  \tilde{I}^n_{\geq m+P}(x) \cap \tilde{I}^m_{\geq m+P}(y)
  =\emptyset.
  \end{equation*}
\end{lem}
This means that, if it not forbidden by $x$ to choose $y$ at time $m$,
then there is no interaction between $x$ and $y$ after time
$m+P$. Thus, the
waiting time $P$ makes it possible to separate completely the two
points (which will be used in Lemma
\ref{interdits_pour_longtemps}). In particular,
  \begin{equation}
  \label{eq_disjointitude}
  I^n_\infty(x)\cap
  I^m_\infty(y)=\emptyset,
  \end{equation}
which implies that the sets we will select in the construction of the
auxiliary partition will be disjoint.
\begin{proof}
Set $U_i=T^n(I^n_\infty(x))$. Assume that $
\tilde{I}^n_{\geq m+P}(x) \cap \tilde{I}^m_{\geq
m+P}(y)\not=\emptyset$, and take a
point $z$ in this intersection. Then
$\dist(T^n z,U_i) \leq \frac{\delta_2}{10}
\sigma^{\frac{m+P-n-1}{2}}$ and $\dist(T^m z, T^m y)\leq
\frac{\delta_2}{10} \left(1+ \sigma^{\frac{P-1}{2}} \right)$.
Note also that, since $y,z\in V_m(y)$, Lemma
\ref{prop_temps_hyperb} implies that $\dist(T^n y, T^n z)\leq
\sigma^{\frac{m-n}{2}} \dist(T^m y, T^m z)$. Hence,
  \begin{align*}
  \dist(T^n y, U_i)
  &
  \leq \dist(T^n y, T^n z) + \dist(T^n z, U_i)
  \leq \sigma^{\frac{m-n}{2}} \dist(T^m y, T^m z) + \dist(T^n z, U_i)
  \\&
  \leq \sigma^{\frac{m-n}{2}} \frac{\delta_2}{10}
  \left(1+\sigma^{\frac{P-1}{2}} \right)
  +\frac{\delta_2}{10} \sigma^{\frac{m+P-n-1}{2}}
  =\frac{\delta_2}{10} \sigma^{\frac{m-n}{2}}
  \left(1+2\sigma^{\frac{P-1}{2}} \right).
  \end{align*}
If $P$ is large enough so that
$1+2\sigma^{\frac{P-1}{2}} \leq
\sigma^{-1/2}$, we get $\dist(T^n y,U_i) \leq
\frac{\delta_2}{10} \sigma^{\frac{m-n-1}{2}}$. As $y\in V_n(x)$
by Lemma \ref{tout_dans_Vn}, we finally get $y\in
\tilde{I}^n_{\geq m}(x)$.
\end{proof}

\begin{lem}
\label{mesure_minoree}
There exists a positive sequence  $c_n$ such that, for all $n\in \N^*$,
for every
$x\in H_n$, $\Leb I^n_\infty(x) \geq c_n$.
\end{lem}
\begin{proof}
The condition $x\in H_n$ implies that, for $k\leq n$, $\dist(T^k
x, S) \geq \alpha_n > 0$, and $T$ is a local diffeomorphism on $M
\moins S$ by definition of $S$. As $T$ is $C^1$ on $\{y \tq
\dist(y, S) \geq \alpha_n\}$, there exists a constant $C_n$ which
bounds $\det DT^n(x)$ for $x \in H_n$. Since the volume-distortion
is bounded by $D_2$ on $V_n(x)$, we get that, for any $y \in
V_n(x)$, $|\det DT^n(y)| \leq D_2 C_n$. In particular, $\Leb
I^n_\infty(x) \geq \Leb(T^n(I^n_\infty(x)))/(D_2 C_n)$. But
$T^n(I^n_\infty(x)))$ is one of the $U_i$, whence its measure is
uniformly bounded away from $0$.
\end{proof}

\begin{lem}
\label{distorsion_hyperbolique_bornee}
There exists a positive constant $C>0$ such that, for any
measurable set $A$, for any $n\in \N^*$, $\Leb(H_n \cap T^{-n}(A) )
\leq C \Leb(A)$.
\end{lem}
\begin{proof}
The sets $I^n_\infty(x)$, for $x\in H_n$, cover $H_n$, and are equal
or disjoint. By Lemma \ref{mesure_minoree}, there is a finite number
of them, say $I^n_\infty(x_1),\ldots, I^n_\infty(x_k)$ (where $k$
depends on $n$).

For $1\leq j \leq k$, the distortion is bounded by $D_2$ on
$I^n_\infty(x_j)$, whence
  \begin{equation*}
  \frac{\Leb(I^n_\infty(x_j) \cap T^{-n}A)}{\Leb(I^n_\infty(x_j))}
  \leq D_2 \frac{\Leb (A)}{\Leb(T^n (I^n_\infty(x_j)))}.
  \end{equation*}
But $T^n(I^n_\infty(x_j))$ is one of the $U_i$, and its measure is
consequently $\geq c$ for some positive $c$. Summing over $j$, we get
  \begin{equation*}
  \Leb(H_n \cap T^{-n}(A)) \leq \frac{D_2}{c} \Leb(A) \Leb(M).
  \qedhere
  \end{equation*}
\end{proof}

\section{The auxiliary partition}

\label{section_Ui}

In this section, we will show the following result (without any
transitivity assumption on $T$):
\begin{thm}
\label{thm_partition_auxiliaire}
Let $T$ be a map on a compact manifold $M$ and
$\delta:\R_+^*\to \R_+^*$ be such that \eqref{non_unif_hyp} and
\eqref{conv_lente} are satisfied. Let $\epsilon^0$ be given by Lemmas
\ref{existe_epsilon_3} and \ref{existe_epsilon_1_2}. We assume that
$T$ satisfies one of the following conditions:
\begin{enumerate}
\item $\Leb \{x \tq h^1_{\epsilon^0}(x)>n\}
=O(u_n)$ where $u_n$ has polynomial decay and tends to $0$.
\item $\Leb\{ x \tq h^1_{\epsilon^0}(x)>n\}=O(u_n)$ where $u_n=e^{-c
n^\eta}$ with $\eta \in (0,1]$.
\item
$\Leb\{ x \tq h^2_{\epsilon^0}(x)>n\}=O(u_n)$
where $u_n$ has polynomial decay and  $(\log n)u_n \to 0$.
\item
$\Leb\{ x \tq h^2_{\epsilon^0}(x)>n\}=O(u_n)$ where $u_n=e^{-c
n^\eta}$ with $\eta \in (0,1]$.
\end{enumerate}
Then there exist a finite partition $U_1,\ldots,U_N$ of $M$, another
finer partition (modulo a set of zero Lebesgue measure)
$W_1,W_2,\ldots$
and times $R_1,R_2,\ldots$ such that, for all $j$,
\begin{enumerate}
\item $T^{R_j}$ is a diffeomorphism between $W_j$ and one of the
$U_i$.
\item $T^{R_j}_{|W_j}$ expands the distances of at least
$\sigma^{-1/2}>1$.
\item The volume-distortion of $T^{R_j}_{|W_j}$ is Lipschitz,
i.e.\ there exists a constant $C$ (independent of $j$) such that, for
every $x,y\in W_j$,
  \begin{equation*}
  \left|1-\frac{\det DT^{R_j}(x)}{\det DT^{R_j}(y)}\right| \leq C
\dist(T^{R_j}x, T^{R_j}y).
  \end{equation*}
\item For $x,y\in W_j$ and $n\leq R_j$, $\dist(T^n x, T^n y) \leq
\dist(T^{R_j}x,T^{R_j}y)$.
\end{enumerate}
Moreover, there exists $c'>0$ such that, under the different
assumptions, the following estimates on the tails hold:
  \begin{equation*}
  \Leb\left( \bigcup_{R_j >n} W_j\right)=\left\{
  \begin{array}{ll}
  O(u_n) & \text{ in the first case},\\
  O((\log n) u_n) &\text{ in the third case},\\
  O(e^{-c' n^\eta}) &\text{ in the second and fourth cases}.
  \end{array}
  \right.
  \end{equation*}
\end{thm}

In the proof of the theorem, it will be sufficient to work on $U_1$,
since the same construction can then be made on each $U_j$.

The fact that the $W_j$ form a partition of $M$ modulo a set of zero
Lebesgue measure will come from the estimates on the size of the
tails, and is not at all trivial from the construction.

This theorem implies the following result on invariant measures:
\begin{thm}
\label{description_mesure_invariante}
Under the assumptions of Theorem
\ref{thm_partition_auxiliaire},
assume moreover that $\sum u_n<\infty$ in the first case, $\sum
(\log n)u_n<\infty$ in the third case.

Then there exists a finite number of invariant absolutely continuous
ergodic probability measures $\mu_1,\dots,\mu_k$.
Moreover, their basins cover almost all
$M$. Finally, there exist disjoint open subsets $O_1,\ldots,O_k$ such
that $\mu_i$ is equivalent to $\Leb$ on $O_i$ and vanishes on $M
\moins O_i$.
\end{thm}
In particular, if $T$ is topologically transitive on $\Lambda$, there
exists a unique absolutely continuous invariant measure.
\begin{proof}[Proof of Theorem
\ref{description_mesure_invariante}]
We build an extension of $M$, similar to a Young tower except that
the basis will be constituted of the finite number of sets
$U_1,\ldots,U_N$. More precisely, set $X=\{ (x,i) \tq x\in W_j,
i<R_j\}$, and let $\pi:X\to M$ be given by $\pi(x,i)=T^i(x)$. We set, for
$x\in W_j$, $T'(x,i)=(x,i+1)$ if $i+1<R_j$, and
$T'(x,R_j-1)=(T^{R_j}(x),0)$. Thus, $\pi \circ T' =T\circ \pi$.
Let $m$ be the measure on $X$ given by $m(A\times\{i\})=\Leb(A)$
when $A\subset W_j$ and $i<j$, so that $\pi_*(m)$ is equivalent to
Lebesgue measure. The condition on the tails ensures that $m$ is
of finite mass.

On $X$, the map $T'$ is Markov, and the map $T'_Y$ induced by $T$ on
the basis $Y=\{(x,0)\}$ is Markov with a Lipschitz volume-distortion
and the big image property. Classical arguments (\cite[Section
4.7]{aaronson:book}) show that $T'_Y$ admits a finite number of
invariant ergodic absolutely continuous probability measures
$\rho_1,\ldots, \rho_l$. Moreover, each of these measures is
equivalent to $m$ on a union $Y_j$ of some sets $U_i \times \{0\}$
(the $Y_j$ are exactly the transitive subsystems for the map $T'_Y$).
Finally, almost every point of $Y$ lands in one of these $Y_j$ after a
finite number of iterations of $T'_Y$.
Inducing (\cite[Proposition 1.5.7]{aaronson:book}), we get a finite
number of absolutely continuous invariant ergodic measures
$\nu_1,\ldots, \nu_l$, whose
basins cover almost all $X$. The condition on the measure of the tails
ensures that the $\nu_i$ are still of finite mass, whence we can
assume that they are probability measures.

The measures $\pi_*(\nu_i)$ are not necessarily all different. Let
$\mu_1,\ldots,\mu_k$ be these measures without repetition. They are
ergodic, and their basins cover almost all $M$, whence there is no
other absolutely continuous invariant ergodic measure.

Let $\mu=\pi_*(\nu)$ be one of the measures $\mu_j$. Since $\nu$ is
equivalent  to $m$ on some set $U_i \times \{0\}$, $\mu$ is equivalent
to $\Leb$ on $U_i$. We will construct the open set $O(\mu)$ of the
statement of the theorem. Let $\Omega_0$ be the interior of $U_i$ (it
is nonempty by construction). By induction, if $\Omega_n$ is defined
and open, set $\Omega_{n+1}=T(\Omega_n \moins S) \cup \Omega_n$. As
$S$ is closed and $T$ is a local diffeomorphism outside of $S$,
$\Omega_{n+1}$ is still an open set.  Set $O=\bigcup \Omega_n$. As
$\mu$ is invariant, we check by induction that $\mu$ is equivalent to
$\Leb$ on $\Omega_n$, whence on $O$. Let us show that, if
 $A\subset M\moins O$, then $\mu(A)=0$. Otherwise, by ergodicity,
there would exist $n$ such that $\mu(T^{-n}(A)\cap \Omega_0)>0$. As
$\mu(S)=0$ (since $\Leb(S)=0$), we get $\mu(T^{-n}(A)\cap (\Omega_0 \moins
S))>0$, whence  $\mu(T^{-(n-1)}(A) \cap \Omega_1)>0$. By induction,
$\mu (A \cap \Omega_n)>0$, which is a contradiction.
\end{proof}

This result is a first step towards the spectral decomposition of
$T$. It was already known, under weaker assumptions (see in
\cite{alves_bonatti_viana} the remark following Corollary D). We will
get later a complete spectral decomposition: each measure $\mu_i$ has
a finite number of components which are mixing (and even exact) for an
iterate of $T$ (Theorem  \ref{vitesse_decroissance}, which also gives
the speed of decay of correlations).

\subsection{Description of the construction}

\label{subsection_construction_Ui}
To prove Theorem \ref{thm_partition_auxiliaire},
we will build a partition of $U_1$ by sets
$W_1,W_2,\ldots$ such that, for every $n$, there exists a return time
$R_n$ such that $T^{R_n}$ is an isomorphism between $W_n$ and one of
the $U_i$, expanding of at least $\sigma^{-n/2}$ and whose
volume-distortion is $D_1$-Lipschitz. In fact, $W_n$ will be some set
$I^{R_n}_\infty(x)$.
Set $H_n(U_1)=H_n \cap \{y \in U_1 \tq \dist(y,\partial U_1)
\geq \lambda_1^n\}$. Hence, if $x\in H_n(U_1)$, we have $V_n(x)
\subset U_1$ by the fourth assertion of Lemma \ref{prop_temps_hyperb}.

We build in fact points
$x^1_1,\ldots, x^1_{l(1)}$ at time $1$, and
$x^2_1,\ldots,x^2_{l(2)}$ at time $2$, and so on.
They will satisfy the following properties:
\begin{itemize}
\item
$x^n_1,\ldots,x^n_{l(n)}$ belong to $H_n(U_1) \moins
\bigcup_{i<n,j\leq l(i)} \tilde{I}^i_{\geq n}(x^i_j)$, and this set is
covered by
$\bigcup_j I^n_\infty(x^n_j)$.
\item
the sets $I^n_\infty(x^n_j)$ (for $n\in \N^*$ and $1\leq j \leq l(n)$)
are disjoint, and
included in $U_1$.
\end{itemize}
We will take for $W_j$ the sets $I^n_\infty(x^n_i)$, and the
corresponding return time $R_j$ will be $n$.

\begin{proof}[Construction of $x^n_i$]
The construction is by induction on $n$. At time $n$, note that, if
$x,y\in H_n(U_1)$, then $I^n_\infty(x)$ and $I^n_\infty(y)$ are either
disjoint or equal. Hence, there exists a system
$I^n_\infty(x^n_1),\ldots, I^n_\infty(x^n_{l(n)})$ of representatives
of the sets
$I^n_\infty(x)$ for $x\in H_n(U_1) \moins \bigcup_{i<n, j\leq
l(i)} \tilde{I}^i_{\geq n}(x^i_j)$ (and it is finite by
Lemma \ref{mesure_minoree}).

By construction, two sets $I^n_\infty(x^n_i)$ constructed at the same
time are disjoint. Take $m>n$, and
$x^m_k \in H_m(U_1)\moins
\bigcup_{i<m,j\leq l(i)} \tilde{I}^i_{\geq n}(x^i_j)$. Then
$x^m_k \in H_m \moins \tilde{I}^n_{\geq m}(x^n_i) $, whence Lemma
\ref{lemme_disjointitude} ensures that $I^m_\infty(x^m_k)$
is disjoint from $I^n_\infty(x^n_i)$.

Finally, to see that $I^n_\infty(x^n_i) \subset U_1$, we use
the fact that $x^n_i \in H_n(U_1)$, whence $\dist(x^n_i, \partial U_1)
\geq \lambda_1^n$. As $V_n(x^n_i) \subset B(x^n_i,
\lambda_1^n)$, this implies that $I^n_\infty(x^n_i) \subset U_1$.
\end{proof}

The properties of hyperbolic times given in Lemma
\ref{prop_temps_hyperb} imply that the expansion and distortion
requirements of Theorem \ref{thm_partition_auxiliaire} are
satisfied. It only remains to estimate
$\Leb\{ x\tq \exists j, x\in W_j \text{ and }R_j >
n\}$.

\subsection{Measure of points which are forbidden many times}
\label{majoration_mesure}

We will denote by $I_n$ the set of points which are forbidden at the
instant $n$, i.e.\
  \begin{equation*}
  I_n=\bigcup_{i<n, j\leq l(i)} \tilde{I}^i_{\geq n}(x^i_j),
  \end{equation*}
and $I^n$ the set of points which are forbidden by the instant $n$, i.e.\
  \begin{equation*}
  I^n=\bigcup_{j\leq l(n)} \tilde{I}^n_{\geq n+1}(x^n_j).
  \end{equation*}
In particular, $I^n \subset I_{n+1}$. Finally, set
  \begin{equation}
  S_n=\bigcup_{i \leq n, j\leq l(i)} I^i_\infty(x^i_j).
  \end{equation}
This is the set of points which are selected before the instant $n$. In this
paragraph, the word ``time'' will be used only for durations, and
``instant'' will be used otherwise.

In this paragraph, we will prove Lemma \ref{conclut_petite_mesure},
which says that the set of points which are forbidden at $k$
instants without being selected has a measure which decays
exponentially fast. The argument is combinatorial: if a point is
forbidden by few instants, then it will be forbidden for a long
time at many of these instants, and it is easily seen that this
gives a small measure (Lemma \ref{interdits_pour_longtemps}).
Otherwise, the point is forbidden by many instants, and we have to
see that each of these instants enables us to gain a
multiplicative factor $\lambda<1$. We will treat two cases: either
the forbidden sets are included one in each other, whence only a
proportion $<1$ is kept at each step, which concludes (Lemma
\ref{interdits_imbriques}), or the forbidden sets intersect each
other close to their respective boundaries, and we just have to
ensure that these boundaries are small enough (Lemma
\ref{interdits_frontiere}).

We will write $B$ for a set $\tilde{I}^n_{\geq n+1}(x^n_i)$, i.e.\
a ``forbidden ball'' (where $x^n_i$ is one of the points defined
in the construction of paragraph
\ref{subsection_construction_Ui}). Then $t(B)$ will denote the
instant $n$ by which it is forbidden, while the ``core''
$C(B)=I^n_\infty(x^n_i)$ is the inner part of $B$, corresponding
to points which are really selected. If $T^{t(B)}(C(B))=U_i$, then
$T^{t(B)}(B)=\{ x \tq \dist(x,U_i) \leq \frac{\delta_2}{10} \}$,
whence $\diam T^{t(B)}(B) \leq \frac{3 \delta_2}{10} \leq
\frac{\delta_2}{2}$. In all the statements and proofs of this
paragraph, the sets denoted by $B_i$ or $B'_i$ will implicitly be
such forbidden balls. We will define in the following lemmas sets
$Z^1,\ldots,Z^6$ of ``points which are forbidden at many instants'', and
we will see that each of them has an exponentially small measure.

\begin{lem}
\label{interdits_frontiere}
Let $Q\in \N^*$. Set
  \begin{multline*}
  Z^1(k,B_0)=\Biggl\{ x \tq \exists B'_1,B_1,\ldots,B'_r,B_r \text{ with }
  \forall 1\leq i \leq r, t(B_{i-1})\leq t(B'_i) \leq  t(B_i)-Q, B_i \not
  \subset B'_i,
  \\
  \sum_{i=1}^r \left \lfloor \frac{t(B_i)-t(B'_i)}{Q}
  \right\rfloor  \geq k, \text{ and }x \in
  \bigcap_{i=0}^r B_i \cap \bigcap_{i=1}^r B'_i \Biggr\}.
  \end{multline*}
Then there exists a constant $C_3$ (independent of $Q$) such that
for all $k$ and $B_0$, $\Leb(Z^1(k,B_0)) \leq C_3 (C_3
\lambda_2^Q)^k \Leb(C(B_0))$.
\end{lem}
Recall that $\lambda_2$ is a constant satisfying
\eqref{frontiere_exp} and \eqref{mesure_frontiere_U'_i}.
\begin{proof}
Let $C_3$ be such that, for $1\leq i\leq N$, $\Leb\{ x \tq
\dist(x,U_i) \leq \frac{\delta_2}{10} \} \leq \frac{C_3}{D_2}
\Leb(U_i)$, and such that $\frac{C_3^{-1}}{1-C_3^{-1}}
(C_2D_2)^2\leq 1$. We will prove that $C_3$ satisfies the
assertion of the lemma, by induction on $k$.

Take $k=0$. Let $n=t(B_0)$, and $i$ be such that
$T^n(C(B_0))=U_i$. Then $Z^1(0,B_0)=B_0$,
whence $T^n(Z^1(0,B_0))= \{ x \tq \dist(x,U_i) \leq
\frac{\delta_2}{10} \}$. This gives $\Leb(T^n(Z^1(0,B_0)))
\leq \frac{C_3}{D_2} \Leb(T^n(C(B_0)))$. As the distortion of
$T^n$ is bounded by $D_2$, by \eqref{distortion_bornee}, we get
$\Leb(Z^1(0,B_0)) \leq C_3 \Leb(C(B_0))$.

Take now $k\geq 1$. Then, decomposing according to the value of
$B'_1$, we get
  \begin{equation*}
  Z^1(k,B_0) \subset \bigcup_{t=1}^k \bigcup_{B'_1 \cap B_0
  \not=\emptyset}
  \bigcup_{\substack{B_1\cap B'_1\not=\emptyset, B_1 \not\subset B'_1
  \\ \left\lfloor
  \frac{t(B_1)-t(B'_1)}{Q}\right\rfloor  \geq t}} Z^1(k-t, B_1).
  \end{equation*}

Let us show that, if $t(B_1)-t(B'_1)=n$, then $B_1$ is included in
an annulus of size $\sigma^{n/2}$ around $B'_1$. More precisely,
set $p=t(B'_1)$, $U'_i=T^p(B'_1)$, and let us show that
  \begin{equation}
  \label{ininfein}
  T^p(B_1) \subset\left\{
  y \tq \dist(y,\partial U'_i) \leq \frac{\delta_2}{2}\sigma^{n/2}
  \right\}.
  \end{equation}
Note that $B_1$ contains a point of $\partial B'_1$, since it is
connected and intersects $B'_1$ and its complement. Thus,
$T^p(B_1)$ contains a point of $\partial U'_i$. Moreover,
\begin{equation*}
  \diam T^p(B_1) \leq \sigma^{n/2} \diam T^{n+p}(B_1) \leq \sigma^{n/2}
  \frac{\delta_2}{2}.
  \end{equation*}
This shows \eqref{ininfein}. Note that
\eqref{mesure_frontiere_U'_i} gives an upper bound for the measure
of \eqref{ininfein}.

Since the distortion is bounded by $D_2$ at hyperbolic times, and
the cores $C(B_1)$ are disjoint by construction, we get by
\eqref{ininfein} and \eqref{mesure_frontiere_U'_i} that
  \begin{equation}
  \sum_{\substack{B_1\cap B'_1\not=\emptyset, B_1 \not\subset B'_1
  \\ \left\lfloor
  \frac{t(B_1)-t(B'_1)}{Q}\right\rfloor  \geq t}} \Leb( C(B_1))
  \leq C_2 \lambda_2^{Qt} D_2 \Leb(C(B'_1)).
  \end{equation}
Finally, write $q=t(B_0)$. Let $x$ be such that
$C(B_0)=I_\infty^q(x)$. The sets $C(B'_1)$ are pairwise disjoint
by construction, and included in $V_q(x)$ by Lemma
\ref{tout_dans_Vn}. Moreover, $T^q$ is a diffeomorphism on
$V_q(x)$ and its distortion is bounded by $D_2$. Since
$T^q(C(B_0))$ is a set $U_i$ and $T^q(V_q(x))= B(T^q x,
\delta_2)$, we have $\Leb (T^q(V_q(x))) \leq C_2 \Leb(
T^q(C(B_0)))$ by \eqref{grosse_mesure_Ui}. By bounded distortion,
we obtain
  \begin{equation}
  \sum_{B'_1 \cap B_0 \not =\emptyset} \Leb (C(B'_1))
  \leq C_2 D_2 \Leb (C(B_0)).
  \end{equation}
Using the induction assumption, we finally obtain
  \begin{align*}
  \Leb Z^1(k,B_0) &
  \leq \sum_{t=1}^k \sum_{B'_1 \cap B_0 \not=\emptyset}
  \sum_{\substack{B_1\cap B'_1\not=\emptyset, B_1 \not\subset B'_1
  \\ \left\lfloor
  \frac{t(B_1)-t(B'_1)}{Q}\right\rfloor  \geq t}} \Leb Z^1(k-t, B_1)
  \\&
  \leq\sum_{t=1}^k \sum_{B'_1 \cap B_0 \not=\emptyset}
  \sum_{\substack{B_1\cap B'_1\not=\emptyset, B_1 \not\subset B'_1
  \\ \left\lfloor
  \frac{t(B_1)-t(B'_1)}{Q}\right\rfloor  \geq t}} C_3 (C_3
  \lambda_2^Q)^{k-t} \Leb (C(B_1))
  \\&
  \leq\sum_{t=1}^k \sum_{B'_1 \cap B_0 \not=\emptyset}
   C_3 (C_3
  \lambda_2^Q)^{k-t} C_2 \lambda_2^{Qt} D_2 \Leb(C(B'_1))
  \\&
  \leq
  C_3 \lambda_2^{Qk} C_3^k (C_2 D_2)^2 \left(\sum_{t=1}^k C_3^{-t}\right)
  \Leb (C(B_0)).
  \end{align*}
By definition of $C_3$, we have $(C_2 D_2)^2 \left(\sum_{t=1}^k
C_3^{-t}\right) \leq 1$. This concludes the induction.
\end{proof}

\begin{lem}
\label{interdits_imbriques}
Set
  \begin{equation*}
  Z^2_{k,N} =\{ x \tq \exists B_1 \varsupsetneq B_2 \cdots
  \varsupsetneq B_k \text{ with
  }t(B_k) \leq N \text{ and }
  x\in (B_1 \cap \cdots \cap B_k) \moins S_N  \}.
  \end{equation*}
Then there exists a constant $\lambda_3<1$ such that
$\Leb(Z^2_{k,N}) \leq \lambda_3^k \Leb(M)$.
\end{lem}
\begin{proof}
We fix $N$ once and for all in this proof, and we will omit all
indexes $N$. We will show that $\lambda_3=\frac{C_2D_2}{C_2D_2+1}$
satisfies the conclusion of the lemma. Note that, for every $B$,
  \begin{equation}
  \label{rapport_centre}
  \Leb(B) \leq C_2 D_2 \Leb(C(B))
  \end{equation}
by \eqref{grosse_mesure_Ui} and the bounded distortion of
hyperbolic times.

We will write  $\boB_1$ for the sets of balls $B$ with $t(B)\leq
N$ which are not included in any other ball $B'$. Write also
$\boB_2$ for the set of balls $B\not\in \boB_1$ with $t(B)\leq N$
which are included only in balls of $\boB_1$, and so on. We will
say that a ball of $\boB_i$ has rank $i$. Every ball $B$ has
finite rank, since a ball which is constructed at time $n$ has at
most rank $n$.

Set $S'_k=\bigcup_{i=1}^k \bigcup_{B \in \boB_i} C(B)$: these are
the points which are selected in balls of rank at most $k$.
Set
  \begin{equation*}
  Z^3_k=\left(\bigcup_{B\in \boB_k} B\right) \moins S'_k.
  \end{equation*}
Let us show that $Z^2_k \subset Z^3_k$.

Take $x\in Z^2_k$, it is in a set $(B_1 \cap \cdots \cap B_k)
\moins S_N$ with $B_1\varsupsetneq B_2 \cdots \varsupsetneq B_k$
and $t(B_k)\leq N$. In particular, $B_k$ is of rank $r\geq k$. Take
$B'_1\varsupsetneq B'_2 \varsupsetneq \cdots \varsupsetneq
B'_{r-1} \varsupsetneq B'_r$ a sequence with $B'_i \in \boB_i$ and
$B'_r=B_k$. In particular, $x\in B'_k$. Moreover, $S'_k \subset
S_N$. As $x\not \in S_N$, we get $x\not\in S'_k$. This shows that
$x\in Z^3_k$.

Let us estimate $\Leb(Z^3_{k+1})$ using $\Leb(Z^3_k)$. Consider
$B_{k+1}\in \boB_{k+1}$.
Let $B_k$ be a ball of rank $k$ containing $B_{k+1}$. As the cores
of different balls are disjoint, $C(B_{k+1}) \cap S'_k=\emptyset$.
Thus, $C(B_{k+1}) \subset B_k \moins S'_k \subset Z^3_k$. However,
$C(B_{k+1}) \subset S'_{k+1}$ by definition, whence $C(B_{k+1})
\cap Z^3_{k+1}=\emptyset$. This shows that $C(B_{k+1}) \subset
Z^3_k \moins Z^3_{k+1}$.

Finally, by \eqref{rapport_centre},
  \begin{equation*}
  \Leb(Z^3_{k+1}) \leq \sum_{B_{k+1} \in \boB_{k+1}} \Leb(B_{k+1})
  \leq C_2D_2 \sum_{B_{k+1} \in \boB_{k+1}} \Leb(C(B_{k+1}))
  \leq C_2D_2 \Leb(Z^3_k \moins Z^3_{k+1})
  \end{equation*}
since the $C(B_{k+1})$ are disjoint. Hence,
  \begin{equation*}
  (C_2 D_2+1) \Leb(Z^3_{k+1}) \leq C_2D_2 \Leb(Z^3_{k+1}) + C_2D_2
  \Leb(Z^3_k \moins Z^3_{k+1}) = C_2D_2 \Leb(Z^3_k).
  \end{equation*}

We obtain by induction that $\Leb(Z^3_k) \leq
\left(\frac{C_2D_2}{C_2D_2+1} \right)^k \Leb(M)$, which gives the
same inequality for $\Leb(Z^2_k)$ since $Z^2_k \subset Z^3_k$.
\end{proof}

\begin{lem}
\label{majore_beaucoup_instants}
Set
  \begin{equation*}
  Z^4(k,N)=\{ x \tq \exists t_1<\ldots<t_k \leq N, x\in I^{t_1}\cap
\cdots \cap I^{t_k} \} \moins S_N.
  \end{equation*}
There exist constants $C_4>0$ and $\lambda_4<1$ such that, for all
$1\leq k\leq N$,
$\Leb(Z^4(k,N)) \leq C_4 \lambda_4^k$.
\end{lem}
This lemma means that the points forbidden by at least $k$ instants
have an exponentially small measure.
\begin{proof}
Take $Q$ large enough so that $C_3 \lambda_2^Q <1$ in Lemma
\ref{interdits_frontiere}. Write $N=rQ+s$ with $s<Q$.

Let $x\in Z^4(k,N)$, forbidden by the instants $t_1<\ldots<t_k$.
For $0\leq u<r$, we choose in each interval $[uQ, (u+1)Q)$ the
first instant $t_i$ (if there exists one), which gives a sequence
$t'_1<\ldots<t'_{k'}$, with $Qk'+s \geq k$. Then we keep the
instants with an odd index, which gives a sequence of instants
$u_1<\ldots<u_l$ with $2l \geq k'$, whence $l \geq k/(2Q)-s$.
Moreover, $u_{i+1}-u_i \geq Q$ for all $i$. Let $B_1,\ldots,B_l$
be balls constructed at the instants $u_i$ and forbidding $x$.

Set $I=\{ 1\leq i \leq l, B_i \subset B_1 \cap \ldots \cap
B_{i-1}\}$ and $J=[1,l]\moins I$. If $\Card I \geq l/2$, we keep only
the balls whose indexes are in $I$. Since there are at least $l/2$
such balls,  $x \in Z^2_{l/2, N}$ (where $Z^2$ is defined in Lemma
\ref{interdits_imbriques}). This lemma implies that the points
obtained in this way have an exponentially small measure (in $l$,
whence in $k$).

Otherwise, $\Card J \geq l/2$. Let $j_0=\sup J$, and $i_0=\inf\{
i<j_0, B_{j_0}\not\subset B_i\}$. Let $j_1=\sup\{ j\leq i_0, j\in
J\}$, and $i_1=\inf\{i<j_1, B_{j_1}\not\subset B_i\}$, and so on: the
construction stops at some step, say $i_n$. Then $J
\subset \bigcup (i_s,j_s]$ by construction, whence $\sum (j_s-i_s)
\geq \Card J \geq l/2$, which implies that $\sum \left \lfloor
\frac{t(B_{j_s})-t(B_{i_s})}{Q} \right\rfloor =\sum \left \lfloor
\frac{u_{j_s}-u_{i_s}}{Q} \right\rfloor \geq l/2$, since two instants
$u_j$ and $u_i$ are separated by at least $Q(j-i)$ by construction.
Hence, the sequence $B_{i_n},B_{i_n}, B_{j_n},\ldots,
B_{i_0},B_{j_0}$ shows that $x \in Z^1(l/2, B_{i_n})$. Summing the
estimates given by Lemma
\ref{interdits_frontiere} over all possible balls $B_{i_n}$,
 we also get an exponentially small measure (since
the cores are disjoint).
\end{proof}

\begin{lem}
\label{interdits_pour_longtemps}
For a ball $B_1=\tilde{I}^{t_1}_{\geq t_1+1}(x_1)$, set
  \begin{multline*}
  Z^5(n_1,\ldots,n_k,B_1)=\{x \tq \exists t_2,\dots,t_k
  \text{ with }t_1<\ldots<t_k\text{ and
  }x_2,\ldots, x_k \text{ such that }\\ \forall 1\leq i\leq k,
  x\in I^{t_i}_{\geq t_i+n_i}(x_i)\}.
  \end{multline*}
There exists a constant $C_5$ (independent of
$B_1,n_1,\ldots,n_k$) such that, when $n_1,\ldots,n_k > P$
(given by Lemma \ref{lemme_disjointitude}),
  \begin{equation*}
  \Leb(Z^5(n_1,\ldots,n_k,B_1)) \leq C_5(C_5 \lambda_2^{n_1})\cdots (C_5
\lambda_2^{n_k}) \Leb(C(B_1)).
  \end{equation*}
\end{lem}
In fact, $Z^5(n_1,\ldots,n_k,B_1)$ is the set of points which are
forbidden for a time at least $n_1$ by $B_1$, and then for a time at
least $n_2$ by another ball $B_2$, and so on.
\begin{proof}
The proof is by induction on $k$.

Let $x\in Z^5(n_1,\ldots,n_k,B_1)$. There exists by definition a
ball $B_2=\tilde{I}^{t_2}_{\geq t_2+1}(x_2)$, constructed at an
instant $t_2>t_1$, such that $x\in Z^5(n_2,\ldots,n_k,B_2)$. The point
$x_2$ is not forbidden at the instant $t_2$ (otherwise, $x_2$ could
not be selected at the instant $t_2$ according to the
construction of paragraph \ref{subsection_construction_Ui}). Hence, Lemma
\ref{lemme_disjointitude} yields that $\tilde{I}^{t_1}_{\geq
t_2+P}(x_1) \cap \tilde{I}^{t_2}_{\geq t_2+P}(x_2)=\emptyset$. But
$x$ is forbidden by the instant $t_2$ for a time at least $n_2>P$,
whence $x\in \tilde{I}^{t_2}_{\geq t_2+P}(x_2)$. Thus, $x\not \in
\tilde{I}^{t_1}_{\geq t_2+P}(x_1)$. As $x\in \tilde{I}^{t_1}_{\geq
t_1+n_1}(x_1)$, we get $t_1+n_1<t_2+P$, i.e.\ $t_2-t_1 > n_1-P$.

Set $U_i=T^{t_1}(C(B_1))$. The expansion at hyperbolic times gives
  \begin{equation*}
  \diam(T^{t_1}(B_2)) \leq \sigma^{\frac{t_2-t_1}{2}}
  \diam(T^{t_2}(B_2)) \leq \sigma^{\frac{n_1-P}{2}} \frac{\delta_2}{2}.
  \end{equation*}
As $\dist(T^{t_1}(x), \partial U_i) \leq \frac{\delta_2}{10}
\sigma^{\frac{n_1-1}{2}}$ since $x$ if forbidden for a time at least
$n_1$, we have proved that there exists a constant $C_6$ such that
  \begin{equation*}
  T^{t_1}(B_2) \subset \boC:=\left\{ y \tq \dist(y, \partial U_i) \leq C_6
\sigma^{\frac{n_1}{2}}\right\}.
  \end{equation*}

By the induction hypothesis, $\Leb(Z^5(n_2,\ldots,n_k,B_2)) \leq
C_5(C_5 \lambda_2^{n_2})\cdots (C_5 \lambda_2^{n_k}) \Leb C(B_2)$.
As the distortion is bounded, we get
$\Leb(T^{t_1}(Z^5(n_2,\ldots,n_k,B_2))) \leq D_2 C_5(C_5
\lambda_2^{n_2})\cdots (C_5 \lambda_2^{n_k})
\Leb(T^{t_1}(C(B_2)))$.
The sets $C(B_2)$ are disjoint by construction and included in
$V_{t_1}(x_1)$ by Lemma \ref{tout_dans_Vn}.
Since $T^{t_1}$ is injective on $V_{t_1}(x_1)$ by Lemma
\ref{prop_temps_hyperb}, the sets $T^{t_1}(C(B_2))$ are still pairwise
disjoint. Moreover, they are all included in the
annulus $\boC$. Hence,
  \begin{align*}
  \Leb(T^{t_1}(Z^5(n_1,\ldots,n_k,B_1))) &
  \leq \sum_{B_2} \Leb(T^{t_1}(Z^5(n_2,\ldots,n_k,B_2)))
  \\&
  \leq C_5 D_2 (C_5 \lambda_2^{n_2})\cdots  (C_5 \lambda_2^{n_k})\sum_{B_2}
  \Leb(T^{t_1}(C(B_2)))
  \\&
  \leq C_5 D_2 (C_5 \lambda_2^{n_2})\cdots  (C_5 \lambda_2^{n_k})
  \Leb(\boC).
  \end{align*}
By \eqref{mesure_frontiere_U'_i}, there exists $C_7$ such that
$\Leb(\boC) \leq C_7 \lambda_2^{n_1} \Leb(U_i)$. Hence,
  \begin{equation*}
  \Leb(T^{t_1}(Z^5(n_1,\ldots,n_k,B_1)))
  \leq C_5 C_7 D_2 \lambda_2^{n_1}
  (C_5 \lambda_2^{n_2})\cdots  (C_5 \lambda_2^{n_k}) \Leb(U_i).
  \end{equation*}
The distortion of the map $T^{t_1}$ is bounded by $D_2$ on
$B_1$. Since $U_i=T^{t_1}(C(B_1))$, the previous equation implies
  \begin{equation*}
  \Leb(Z^5(n_1,\ldots,n_k,B_1)) \leq C_5 C_7 D_2^2 \lambda_2^{n_1}
  (C_5 \lambda_2^{n_2})\cdots  (C_5 \lambda_2^{n_k}) \Leb(C(B_1)).
  \end{equation*}
This concludes the proof, if $C_5 \geq C_7 D_2^2$ is taken large
enough so that the result holds for $k=0$.
\end{proof}

The following lemma will subsume all the previous lemmas: it shows
that the points forbidden at $k$ instants have an exponentially small
measure.
\begin{lem}
\label{conclut_petite_mesure}
Set
  \begin{equation*}
  Z^6(k,N)=\{x \tq \exists t_1<\ldots<t_k \leq N,x \in I_{t_1}\cap
\ldots \cap I_{t_k} \} \moins S_N.
  \end{equation*}
There exist constants $C_8>0$ and $\lambda_5<1$ such that, for all
$k\leq N$,
  \begin{equation*}
  \Leb(Z^6(k,N)) \leq C_8 \lambda_5^k.
  \end{equation*}
\end{lem}
\begin{proof}
Take $R>P$ (given by Lemma \ref{lemme_disjointitude}) so that
$\lambda_2+C_5 \lambda_2^R < 1$. Let $x\in Z^6(k,N)$, and consider
all the instants $u_i$ by which it is forbidden for a time
$n_i\geq R$, ordered so that $u_1<\dots<u_p$. Then $x\in
Z^5(n_1,\ldots,n_p,B_1)$ for some ball $B_1$. If $\sum n_i \geq
k/2$, we do not do anything else. Otherwise, let $v_1<\ldots<v_q$
be the other instants by which $x$ is forbidden, for times
$m_1,\dots,m_q<R$. Then $\sum n_i+\sum m_j$ is not less than the
number of instants at which $x$ is forbidden, whence $\sum m_j
\geq k/2$. This implies that $Rq\geq k/2$. We obtain
  \begin{equation*}
  Z^6(k,N) \subset \left(\bigcup_{B_1} \bigcup_{\substack{n_1,\ldots,n_p
\geq R \\ \sum n_i \geq k/2}}Z^5(n_1,\ldots,n_p,B_1)\right) \cup
Z^4(k/(2R),N).
  \end{equation*}
Consequently, Lemmas \ref{majore_beaucoup_instants} and
\ref{interdits_pour_longtemps} yield that
  \begin{equation*}
  \Leb (Z^6(k,N)) \leq \sum_{B_1} \sum_{\substack{n_1,\ldots,n_p
  \geq R \\ \sum n_i \geq k/2}}C_5(C_5 \lambda_2^{n_1}) \cdots
  (C_5 \lambda_2^{n_p}) \Leb(C(B_1)) + C_4 \lambda_4^{k/(2R)}.
  \end{equation*}

As the cores $C(B_1)$ are disjoint, $\sum \Leb(C(B_1)) \leq
\Leb(M)<\infty$. To conclude, it is therefore sufficient to prove
that
  \begin{equation*}
  \sum_{\substack{n_1,\ldots,n_p\geq R\\ \sum n_i \geq k/2}}
  (C_5 \lambda_2^{n_1}) \cdots (C_5 \lambda_2^{n_p})
  \end{equation*}
decays exponentially fast.

We use generating series:
  \begin{equation*}
  \sum_n \sum_{\substack{n_1,\ldots,n_p\geq R\\ \sum n_i=n}}
   (C_5 \lambda_2^{n_1}) \cdots (C_5
  \lambda_2^{n_p}) z^n
  =\sum_{p=1}^\infty\left( C_5 \sum_{n=R}^\infty \lambda_2^n z^n \right)^p
  =\frac{C_5 \lambda_2^R  z^R}{1-\lambda_2 z-C_5 \lambda_2^R  z^R}.
  \end{equation*}
As $\lambda_2+C_5 \lambda_2^R<1$, this function has no pole in a
neighborhood of the unit disk in $\C$. Hence, its coefficients decay
exponentially fast, i.e.\ there exist constants $C_9>0$ and
$\lambda_6<1$ such that
  \begin{equation*}
  \sum_{\substack{n_1,\ldots,n_p\geq R\\ \sum n_i=n}}
   (C_5 \lambda_2^{n_1}) \cdots (C_5
  \lambda_2^{n_p})
  \leq C_9 \lambda_6^n.
  \end{equation*}
We just have to sum over $n\geq k/2$ to conclude.
\end{proof}

\subsection{Proof of Theorem \ref{thm_partition_auxiliaire}}
\label{prouve_thm_auxiliaire}

We check in the four cases of Theorem
\ref{thm_partition_auxiliaire} that the conclusions on the measures of
the
tails hold. In this proof, the precise choice of
$\sigma, \delta$ and $\epsilon^0$ in paragraph
\ref{section_frequence_hyperboliques} is important.
From the previous paragraph, we will only use Lemma
\ref{conclut_petite_mesure}.

\begin{proof}[Proof of the first and second cases]
Recall that $\Leb \{x \tq h^1_{\epsilon^0}(x)>n\} =O(u_n)$. Recall
also that $S_n$ is the set of points selected before time $n$, and
that $\theta$ is defined in Lemma
\ref{existe_epsilon_1_2}. Let us show that
  \begin{equation*}
  U_1 \moins S_n \subset \left\{x \in U_1 \tq h^1_{\epsilon^0}(x)>n\right\}
  \cup \left\{ x\in U_1 \tq \dist(x,
  \partial U_1) \leq \lambda_1^{\theta n/2} \right\}
  \cup Z^6(\theta n/2, n).
  \end{equation*}
This will conclude the proof, since the second and third sets have an
exponentially small measure, by
\eqref{frontiere_exp} and Lemma \ref{conclut_petite_mesure}.

Take $x$ in $U_1 \moins S_n$, which does not belong either to
$\{h^1_{\epsilon^0}(x)>n\}$ or to $\left\{ \dist(x, \partial U_1)
\leq \lambda_1^{\theta n/2} \right\}$. By Lemma
\ref{existe_epsilon_1_2}, $x$ has at least $\theta n$ hyperbolic
times between $1$ and $n$, whence at least $\theta n/2$ between
$\theta n/2$ and $n$. We will denote them by $t_1<\ldots<t_k\leq
n$. As $\dist(x,\partial U_1)> \lambda_1^{\theta n/2}$, we have in
fact $x\in H_{t_i}(U_1)$ for all these instants. If $x$ was not
forbidden at the instant $t_i$, then it would be selected at the
instant $t_i$ by construction, which is not possible since
$x\not\in S_n$. Hence, $x\in I_{t_i}$. We obtain in this way at
least $\theta n/2$ instants at which $x$ is forbidden, whence
$x\in Z^6(\theta n/2, n)$.
\end{proof}

\begin{proof}[Proof of the third and fourth case]
Denote by $N(x,n)$ the number of hyperbolic times of $x$ between
$1$ and $n$.
\begin{lem}
\label{lemme_combine}
Let $n\in \N^*$ and $k(n) \in [1, \theta n]$. Then
  \begin{equation*}
  \Leb\{x\tq N(x,n) < k(n)\}
  \leq C \frac{k(n)}{\theta} \Leb\left\{ x \tq h^2_{\epsilon^0}(x)>
  n-\frac{k(n)}{\theta}\right\}.
  \end{equation*}
\end{lem}
\begin{proof}
Write $SH_l^*$ for the set of points whose first positive super hyperbolic
time is $l$. If a point $x$ has a super hyperbolic time $j$
between $k(n)/\theta$ and $n$, then it will have at least $\theta
j \geq k(n)$ hyperbolic times between $1$ and $j$, by Lemma
\ref{existe_epsilon_1_2}. Hence,
  \begin{equation*}
  \{x\tq  N(x,n) < k(n) \} \subset M \moins \bigcup_{k(n)/\theta
  \leq j \leq n} SH_j.
  \end{equation*}
Denote by $k \in [0, k(n)/\theta)$ the last super hyperbolic time
of $x$ before $k(n)/\theta$. We get
  \begin{equation*}
  \Leb\left(M \moins \bigcup_{k(n)/\theta
  \leq j \leq n} SH_j\right)
  \leq \sum_{k=0}^{k(n)/\theta}
  \Leb\left( SH_k\cap T^{-k} \left(\bigcup_{l>n-k}
  SH_l^*\right) \right)
  \leq C \sum_{k=0}^{k(n)/\theta}
  \Leb \left(\bigcup_{l>n-k}
  SH_l^*\right),
  \end{equation*}
using the inclusion $SH_k \subset H_k$ and Lemma
\ref{distorsion_hyperbolique_bornee} for the last inequality.

By Lemma \ref{existe_epsilon_3}, a point $x$ has at least one
super hyperbolic time between $1$ and $h^2_{\epsilon^0}(x)$,
whence $\bigcup_{l>n-k} SH_l^* \subset \{x\tq
h^2_{\epsilon^0}(x)>n-k\}$. This concludes the proof of the lemma.
\end{proof}
For any $k(n)$, the same arguments as in the proof of the first
and second cases imply that
  \begin{equation*}
  \label{inclusion_U1Sn}
  U_1 \moins S_n \subset \{ x \tq N(x,n) < k(n)\} \cup
  \{ x \tq \dist(x,\partial U_1) \leq \lambda_1^{k(n)/2}\}
  \cup Z^6(k(n)/2,n).
  \end{equation*}
By \eqref{frontiere_exp}, Lemma \ref{conclut_petite_mesure} and
Lemma \ref{lemme_combine}, we get
  \begin{equation}
  \label{lqsjfdlkjsqdflk}
  \Leb(U_1 \moins S_n) \leq C \frac{k(n)}{\theta}
  \Leb \left\{ h^2_{\epsilon^0}(x) >
  n-\frac{k(n)}{\theta} \right\} + C_2 \lambda_2^{k(n)/2}+C_8
  \lambda_5^{k(n)/2}.
  \end{equation}
To conclude the proof, we just have to choose correctly the
sequence $k(n)$.

Assume that $\Leb \{x \tq h^2_{\epsilon^0}(x)>n\} =O(u_n)$ where
$u_n$ has polynomial decay. Choose $K$ large enough so that
$k(n):= \lfloor K \log n \rfloor$ satisfies $\lambda_5^{k(n)/2}=O(u_n)$ and
$\lambda_2^{k(n)/2}=O(u_n)$. Then \eqref{lqsjfdlkjsqdflk} gives
$\Leb(U_1 \moins S_n) =O ( (\log n)u_{n-k(n)/\theta}) = O((\log
n)u_n)$.

Assume finally that $\Leb \{x \tq h^2_{\epsilon^0}(x)>n\} =O(e^{-c
n^\eta})$ with $\eta\in (0,1]$. Choose $k(n)=\lfloor n^\eta \rfloor$
if $\eta<1$,
and $k(n)= \left \lfloor \frac{\theta}{2}n \right\rfloor$ if $\eta=1$. Then
\eqref{lqsjfdlkjsqdflk} gives $\Leb(U_1 \moins S_n)=O(e^{-c'
n^\eta})$ for some $c'>0$.
\end{proof}
The logarithmic loss in the polynomial case comes from the factor
$k(n)$ in Lemma \ref{lemme_combine}.

\section{The Young tower}
\label{section_tour}

Using Theorem \ref{thm_partition_auxiliaire}, it is possible to
prove directly the estimates on the decay of correlations (under a
mixing assumption): the coupling arguments of
\cite{lsyoung:recurrence} apply to the ``tower'' built from the
partition $W_j$ (the only difference with the towers of
\cite{lsyoung:recurrence} is that the returns to the basis do not
cover the whole basis, but only one of the sets $U_i$). This is
for example shown in \cite{gouezel:these}. However, in view of the
existing literature, it seems more economical to build a true
Young tower, in order to apply directly the results of
\cite{lsyoung:recurrence} (or rather a small improvement of these
results, since the results of Young are not sharp enough in the
stretched exponential case).

\subsection{Construction of the Young tower}

The Young tower is given by the following theorem:
\begin{thm}
\label{thm_partition_young}
Under the assumptions of Theorem
\ref{description_mesure_invariante}, let $\mu$ be one of the invariant
absolutely continuous ergodic probability measures given by this theorem.

Then there exist a nonempty open set $B$ on which $\mu$ is equivalent
to Lebesgue measure, a partition (modulo $0$)
$Z_1,Z_2,\ldots$ of $B$, and times $R'_1,R'_2,\ldots$ such that, for all
$j$
\begin{enumerate}
\item $T^{R'_j}$ is a diffeomorphism between $Z_j$ and $B$.
\item $T^{R'_j}_{|Z_j}$ expands the distances of at least
$\sigma^{-1/2}>1$.
\item the volume-distortion of $T^{R'_j}_{|Z_j}$ is Lipschitz.
\item For $x,y\in Z_j$ and $n\leq R'_j$, $\dist(T^n x, T^n y) \leq
\dist(T^{R'_j}x,T^{R'_j}y)$.
\end{enumerate}
Moreover, the estimates on the size of tails as given in Theorem
\ref{thm_partition_auxiliaire} still hold.
\end{thm}

\begin{proof}
Let $X$ be the extension of $M$ constructed in the proof of Theorem
\ref{description_mesure_invariante} using the auxiliary partition, and
$\nu$ one of the invariant ergodic measures on $X$ such that
$\pi_*(\nu)=\mu$. We identify each set $U_i$ in $M$
with $U_i\times \{0\}$ in $X$.

On one  $U_i$ (let us say $U_1$), the measure $\nu$ is equivalent to
$m$. The basis $B$ of the Young tower will be $U_1$. Write
$U_2,\ldots,U_s$ for the other sets $U_i$ on which $\nu$ is equivalent
to $m$.
Let $T'_Y$ be the map induced by $T'$ on
$Y=\{(x,0)\} \subset X$, i.e., on an element $W_j$ of the partition
$\boB$ given by Theorem \ref{thm_partition_auxiliaire}, with return time
$R_j$, we set $T'_Y(x,0)=(T^{R_j}(x),0)$. We define a partition
$\boB^n$ of $Y$ by $\boB^n =\bigcap_0^{n-1} (T'_Y)^{-i}
(\boB)$: thus, an element of $\boB^n$ is sent by
$T'_Y,\ldots,(T'_Y)^{n-1}$ on subsets of elements of
$\boB$,
and by $(T'_Y)^n$ on a set $U_i$. As $\nu$ is ergodic,
there exists $L>0$ such that every $U_i$ (with $i\leq s$) contains an
element of $\boB^n$, for some $n<L$, whose image under
$(T'_Y)^n$ is $U_1$.

For $x \in \bigcup_1^s U_i$, we define a sequence of times
$t_0(x)=0,t_1(x), t_2(x), \ldots$ and an integer $k(x)$
(corresponding to the number of iterations before $x$ is selected)
in the following way: let $B_0\in \boB$ contain $x$, and let $R_1$
be its return time. Set $t_1(x)=R_1$. If $T'_Y(B_0)=U_1$, we set
$k(x)=1$ and we stop here. Otherwise, $T^{R_1}(B_0)$ is one of the
sets $U_i$ with $2\leq i \leq s$. We consider the set $B_1$ of the
partition $\boB$ containing $T^{R_1}(x)$, with a return time
$R_2$. Set $t_2(x)=t_1(x)+R_2$. If $T^{t_2(x)}(x)$ is in $U_1$, we
set $k(x)=2$ and we stop here. Otherwise we consider the next
iterate of $T^{t_2(x)}$, that we denote by $T^{t_3(x)}$, and we go
on. More formally, $k(x)=k(T^{R_1}x)+1$ and
$t_j(x)=t_{j-1}(T^{R_1}(x))+t_1(x)$ for every $j\leq k(x)$. By
definition, $k(x)$ is the smallest integer $n\geq 1$ such that the
element of $\boB^n$ containing $x$ is sent on  $U_1$ by
$(T'_Y)^n$.

The elements of the final partition will be the sets $Z_j$
constructed in this way, included in $U_1$,
and the corresponding return time will be
$t_{k(x)}(x)$ for $x\in Z_j$ (this is independent of $x$). By
construction, $T^{t_{k(x)}(x)}(Z_j)=U_1$, and we have a Young
tower.

In the end, almost every point will be selected (we will see later
that the measure of the tails tends to $0$).
The distortion and expansion properties of the
partition $\boB$ ensure that these properties will remain
satisfied by the Young tower. We just have to prove the estimates
on the measures of the tails to conclude.

Set $\tau(x)=t_{k(x)}(x)$. In at most $L$ steps, an element of every
$U_i$ is selected to come back to $U_1$, by definition of $L$.
Since the distortion is bounded, there exists $\epsilon>0$ such
that
\begin{equation}
  \label{young1}
  \Leb(\tau=t_j \text{ or }\ldots \text{ or }\tau=t_{j+L-1}
   \tq t_1,\ldots,t_{j-1}, \tau>t_{j-1}) \geq
  \epsilon.
  \end{equation}
Moreover, still by bounded distortion,
  \begin{equation}
  \label{young2}
  \Leb\{ t_{j+1}-t_j>n \tq t_1,\ldots,t_j\}
  \leq C \sum_{W_k \in \boB, R_k >n } \Leb(W_k),
  \end{equation}
this last term being estimated by Theorem
\ref{thm_partition_auxiliaire}. We want to obtain estimates on the
measure of the tails, i.e.\ on $\Leb\{x\tq \tau(x)>n\}$, and we
will use \eqref{young1} and \eqref{young2} to get them. The
following lemma is indeed sufficient to conclude the proof.
\end{proof}

\begin{lem}
\label{lemme_ameliore_Young}
Let $(X,\mu)$ be a space endowed with a finite measure and $k:X
\to \N$ and $t_0,t_1, t_2,\ldots : X\to \N$ measurable functions such
that $0=t_0<t_1<t_2<\dots$ almost everywhere. Set $\tau(x)=t_{k(x)}(x)$,
and assume that there exist $L>0$ and $\epsilon>0$ such that
  \begin{equation}
  \label{young1prime}
  \mu\{\tau=t_j \text{ or }\ldots \text{ or }\tau=t_{j+L-1}
   \tq t_1,\ldots,t_{j-1}, \tau>t_{j-1}\} \geq
  \epsilon.
  \end{equation}
Assume moreover that there exist a positive sequence $u_n$ and a
constant $C$ such that
  \begin{equation}
  \label{young1seconde}
  \mu\{ t_{j+1}-t_j>n \tq t_1,\ldots,t_j\}
  \leq C u_n.
  \end{equation}
Then
  \begin{enumerate}
  \item  If $u_n$ has polynomial decay, $\mu\{ \tau > n\}=O(u_n)$.
  \item If $u_n=e^{-cn^\eta}$ with $c>0$ and $\eta\in (0,1]$, then
  there exists $c'>0$ such that $\mu\{ \tau >
  n\}=O(e^{-c'n^\eta})$.
  \end{enumerate}
\end{lem}
\begin{proof}
In \cite{lsyoung:recurrence}, Young considers a problem which is
\emph{a priori} completely different: she wants to estimate the
speed of decay of correlations in towers. However, she introduces
a sequence of times $t_n(x)$ which satisfies the assumptions of
the lemma, and she uses only the properties \eqref{young1prime}
and \eqref{young1seconde} to obtain estimates on the set $\mu\{
\tau > n\}$. In particular, in the fourth section of
\cite{lsyoung:recurrence}, she proves our lemma when
$u_n=e^{-cn}$, and when $u_n$ has polynomial decay. She assumes
$L=1$, but her proofs can easily be adapted to the general case.
Moreover, for the polynomial case, she only deals with the case
$u_n=1/n^\gamma$, but the same proof works directly in the general
case, using that $u_{n/i} \leq u_n i^\gamma$  for some $\gamma>0$.

However, in the stretched exponential case (i.e.\ $0<\eta<1$), the
estimates of Young give only $\mu\{\tau >n\}=O(e^{-n^{\eta'}})$
for any $\eta'<\eta$, which is weaker than the result of our
lemma. We will give a different proof in this case.

When $w^1$ and $w^2$ are two real sequences, we will write
$w^1\star
w^2$ for their convolution, given by $(w^1 \star w^2)_n =\sum_{a+b=n}
w^1_a w^2_b$. When $w$ is a sequence, we will also write $w^{\star l}$
for the sequence obtained by convolving $l$ times $w$ with itself.

Write $v_n=C e^{-cn^\eta}$, so that $\mu\{ t_j-t_{j-1} =n \tq
t_{j-1},\dots,t_1\} \leq v_n$.
 Let us show that, for large enough
$K$, the sequence $w_n=1_{n\geq K} v_n$ satisfies
  \begin{equation}
  \label{jsdkjsf}
  \forall p\in \N, (w\star w)_p \leq w_p.
  \end{equation}
Note that, on $[0,1/2]$, the function $(x^\eta+(1-x)^\eta
-1)/x^\eta$ is continuous (it tends to $1$ at $0$), and
positive, whence larger than some constant $\gamma>0$.
Hence, $x^\eta+(1-x)^\eta \geq 1+\gamma x^\eta$. For $p< 2K$,
$(w\star w)_p=0$. Take $p \geq 2K$. Then
  \begin{equation*}
  (w \star w)_p
  \leq 2 C^2 \sum_{K \leq j \leq p/2}e^{-c j^\eta}e^{-c (p-j)^\eta}
  = 2 C^2 \sum_{K\leq j \leq p/2}
      e^{-c p^\eta( (j/p)^\eta+(1-j/p)^\eta ) }.
  \end{equation*}
For $x=j/p$, we have $x\in [0,1/2]$, whence
  \begin{equation*}
  (w \star w)_p
  \leq  2 C^2 \sum_{K\leq j \leq p/2}
      e^{-c p^\eta(1+\gamma (j/p)^\eta) }
  \leq 2 C^2 e^{-c p^\eta} \sum_{j \geq K}  e^{-c \gamma j^\eta}.
  \end{equation*}
Taking $K$ large enough so that $2 C \sum_{j \geq K}  e^{-c \gamma
j^\eta} \leq 1$, we obtain \eqref{jsdkjsf}.

Let $k\geq 0$ and $A \subset \{1,\ldots,k\}$. For $j\in A$, take
$n_j \geq 1$. Set $Y(A,n_j)=\{x \tq k(x) \geq \sup A \text{ and }
\forall j \in A, t_j(x)-t_{j-1}(x)=n_j\}$. Conditioning
successively with respect to the different times, we get
  \begin{equation*}
  \mu \bigl( Y(A,n_j)\bigr) \leq \prod_{j\in A} \mu\{ t_j-t_{j-1}=n_j \tq
  t_{j-1},\ldots,t_1\}
  \leq \prod_{j\in A} v_{n_j}
  \end{equation*}
by \eqref{young1seconde} and the definition of $v_n$.

Set $q(n)=\left\lfloor\alpha n^\eta\right\rfloor$, where $\alpha$ will
be chosen later. Take $x$ such that $\tau(x)>n$. If $k(x)>q(n)$, i.e.\
$x$
is selected after more than $q(n)$ steps, we do not do
anything. Otherwise, let
$l=k(x) \leq
q(n)$, and let $n_j=t_j(x)-t_{j-1}(x)$ for $j \leq l$. Write
$A=\{ j \tq n_j \geq K\}$. Thus, $x \in Y(A,n_j)$. Moreover, as
$\sum n_j=\tau(x)>n$, we have $\sum_{j\in A} n_j \geq n-K q(n)\geq
n/2$ if $n$ is large enough. We have shown that
  \begin{equation}
  \label{T_pas_gros}
  \{ x \tq \tau(x) > n\}
  \subset \{ k(x)> q(n)\} \cup
  \bigcup_{A \subset \{1,\ldots,q(n)\}}
  \bigcup_{\substack{n_j \geq K \\ \sum_A n_j \geq n/2}} Y(A, n_j).
  \end{equation}
By \eqref{young1prime}, $\mu \{ k(x)> q(n)\} \leq
(1-\epsilon)^{q(n)/L} \leq e^{-c'' n^\eta}$ for some $c''$.
Moreover, writing $l=\Card A$ and using \eqref{jsdkjsf},
  \begin{align*}
  \mu \left(
  \bigcup_{A \subset \{1,\ldots,q(n)\}}
  \bigcup_{\substack{n_j \geq K \\ \sum_A n_j \geq n/2}} Y(A, n_j) \right)
  &
  \leq \sum_{A \subset \{1,\ldots,q(n)\}} \sum_{\substack{n_j \geq K
  \\ \sum_A n_j \geq n/2}}  \prod_{j\in A} v_{n_j}
  \\&
  \leq  \sum_{0 \leq l \leq q(n)} \binom{q(n)}{l} \sum_{\substack{
  n_1,\ldots,n_l\geq K \\ \sum n_j \geq n/2}}  v_{n_1} \cdots v_{n_l}
  \\&
  =  \sum_{0 \leq l \leq q(n)} \binom{q(n)}{l}
  \sum_{n/2}^\infty \left(w^{\star l} \right)_p
  \\&
  \leq  \sum_{0 \leq l \leq q(n)} \binom{q(n)}{l}
  \sum_{n/2}^\infty w_p
  = 2^{q(n)} \sum_{n/2}^\infty w_p.
  \end{align*}
As $w_n=O(e^{-c n^\eta})$, one proves (comparing to an integral)
that $\sum_{n/2}^\infty w_p =O(n^{1-\eta} e^{-c (n/2)^\eta})$.
Hence, if $\alpha$ is small enough, $2^{q(n)} \sum_{n/2}^\infty
w_p =O(e^{-c' n^\eta})$ for some $c'>0$. By \eqref{T_pas_gros}, we
have proved that $\mu\{\tau(x)>n\}=O(e^{-c' n^\eta})$.
\end{proof}

\subsection{Consequences}
\label{section_decorrelation}

\begin{thm}
\label{vitesse_decroissance}
Let $T$ satisfy the assumptions of Theorem
\ref{description_mesure_invariante}, $\mu$ be one of the invariant
ergodic absolutely continuous probability measures given by this
theorem, and $O$ be an open set such that $\mu$ is equivalent to
$\Leb_{|O}$.

Then there exists a finite partition (modulo $0$)
$\Omega_0,\ldots,\Omega_{d-1}$
of $O$ in open sets, such that
$T(\Omega_i)=\Omega_{i+1}$ (modulo $0$)
for $i\leq d-1$ ($\Omega_d$ is identified with $\Omega_0$), and such
that,
on each $\Omega_i$, the map $T^d$ is mixing (and even exact) for the
measure $\mu$.

Finally, for every functions $f,g:M\to \R$ with $f$ H\"older and $g$
bounded, there exists a
constant $C$ such that, for $0\leq i\leq d-1$, for all $n\in \N$, the
correlations
$\Cor_{\Omega_i}(f,g\circ
T^{dn}):=\int_{\Omega_i} f\cdot g\circ T^{dn}\dd\mu
-\left(\int_{\Omega_i}f\dd \mu \right) \left(\int_{\Omega_i} g
\dd\mu\right)$ satisfy
  \begin{equation}
  \label{decroit_cor}
  \left|\Cor_{\Omega_i}(f,g\circ T^{dn})\right| \leq \left\{
  \begin{array}{ll}
  C\sum_{p=n}^\infty u_p & \text{ in the first case,}
  \\
  C \sum_{p=n}^\infty (\log p)u_p &\text{ in the third case,}
  \\
  Ce^{-c'n^{\eta}} &\text{ in the second and fourth cases}.
  \end{array}\right.
  \end{equation}
\end{thm}
When all the iterates of $T$ are topologically transitive, there exist a
unique measure $\mu$ and a unique set $\Omega$. This proves Theorems
\ref{thm_principal_0}
et \ref{thm_principal}.

\begin{proof}
Theorem \ref{thm_partition_young} makes it possible to construct an
abstract Young tower
$X=\{(x,i) \tq x\in Z_j, i<R'_j\}$, a
projection $\pi:X\to M$ given by $\pi(x,i)=T^i(x)$, and
a map $T'$ on $X$ such that $\pi\circ T'=T\circ \pi$, as in the proof
of Theorem \ref{description_mesure_invariante} (but using the
partition given by Theorem
\ref{thm_partition_young} instead of the partition given by Theorem
\ref{thm_partition_auxiliaire}).

By \cite{lsyoung:recurrence}, $T'$ admits a unique absolutely
continuous invariant probability measure $\nu$. The measure
$\pi_*(\nu)$ is absolutely continuous with respect to $\mu$, whence
$\pi_*(\nu)=\mu$ by ergodicity.

Set $d_1=\gcd(R'_j)$, and write, for $0\leq k\leq d_1-1$,
$X_k=\{(x,i)\in X \tq i\equiv k \mod d_1\}$. Thus, $T'$ maps
$X_k$ to $X_{k+1}$ for $k<d_1$ (taking $k$ modulo $d_1$).
The system $(X_k, (T')^{d_1})$ is then a Young tower whose return
times are relatively prime, and whose invariant measure is
$\nu_k:=\nu_{|X_k}$. \cite[Theorem
1]{lsyoung:recurrence} implies that $\nu_k$ is exact for
$(T')^{d_1}$. Moreover, the correlations of H\"older functions (as
defined in \cite{lsyoung:recurrence}) decay as indicated in
\eqref{decroit_cor}: in the exponential case, this is proved in
\cite{lsyoung:recurrence}. Young treats the case of $1/n^\gamma$, but
her proof can easily be adapted to the polynomial case. It remains to
treat the stretched exponential case, which is given by the following lemma:
\begin{lem}
Let $(X,T')$ be a mixing Young tower, and assume that the return
time on the basis $R$ satisfies $m(R>n)=O(e^{-c n^\eta})$ for some
$0<\eta<1$. Then, if $f$ is H\"older and $g$ is bounded, the
correlations of $f$ and $g$ are bounded by $e^{-c' n^\eta}$ for
some $c'>0$.
\end{lem}
\begin{proof}
This is a consequence of \cite[Section 3.5]{lsyoung:recurrence} and Lemma
\ref{lemme_ameliore_Young}.
\end{proof}

These results are true on $X$, we still have to come back to $M$.

The measures $\lambda_k=\pi_*(\nu_k)$ satisfy $T_*
\lambda_k=\lambda_{k+1}$, and are invariant and ergodic for
$T^{d_1}$. In particular, two such measures are either equal or
mutually singular. Hence, there exists $d$ (dividing $d_1$, let us say
$d_1=sd$)
such that $\lambda_k=\lambda_l$ if and only if $k\equiv
l \mod d$. Using the same argument as in the proof of Theorem
\ref{description_mesure_invariante}, we check that the measures
$\lambda_k$ (for $0\leq k<d$) are supported on disjoint open sets
$\Omega_k$. Moreover,
$T_*(\lambda_k)=\lambda_{k+1}$, whence $T(\Omega_k)=\Omega_{k+1}$ modulo $0$.

Let us show that $\lambda_k$ is exact for $T^d$. Let $A\subset
\Omega_k$ have nonzero measure, such that $A$ can we written as
$T^{-dn}(A_n)$ for any $n$. Hence, $A'=\pi^{-1}(A)$ is equal to
$(T')^{-dn}(A'_n)$, where $A'_n=\pi^{-1}(A_n)$. In particular, since
$X_k$ is invariant under $(T')^{d_1}$, we get $A'\cap
X_k=(T')^{-nd_1}(A'_{s n}\cap X_k)$. As $(X_k,\nu_k)$ is
exact, this proves that $A'\cap X_k$ has full $\nu_k$-measure, which
concludes the proof.

Let finally $f,g$ be two functions on $M$ such that $f$ is H\"older
and $g$ is bounded. Write
$f'=f\circ \pi$ and $g'=g\circ \pi$: the function $f'$ is H\"older on
$X$, and $g'$ is bounded. For $n\in \N$, write $n=ps+r$ with
$0\leq r<s$. Then
  \begin{equation*}
  \int_{\Omega_k} f\cdot g\circ T^{dn}=\int_{X_k} f'\cdot (g'\circ
  (T')^{dr})\circ (T')^{p ds}
  =\int_{X_k}f'\cdot (g'\circ
  (T')^{dr})\circ (T')^{p d_1}.
  \end{equation*}
The function $g'\circ (T')^{dr}$ is bounded on $X_k$, whence
the
estimate on the speed of decay of correlations for $\nu_k$ on $X_k$
gives the same estimate for the decay of correlations of $f$ and $g$ on $M$.
\end{proof}

\bibliography{biblio}
\bibliographystyle{alpha}

\end{document}